\documentclass[11pt,leqeq]{article}
\usepackage{amssymb,amsthm}
\usepackage{amsmath}
\usepackage{amscd}
\usepackage{xy}
\xyoption{all}
\usepackage[dvips]{graphicx}
\newtheorem{thm}{Theorem}
\newtheorem{lem}[thm]{Lemma}
\newtheorem{cor}[thm]{Corolary}
\newtheorem{prop}[thm]{Proposition} 

\newtheorem{defn}[thm]{Definition}
\newtheorem{exmp}[thm]{Example}
\newcommand{\des}{\displaystyle}
\newcommand{\id}{{\rm{id}}}

\newcommand{\ob}{{\rm{Ob}}}

\newcommand{\adn}{{\rm{ADN}}}
\newcommand{\sig}{{\rm{sgn}}}

\newcommand{\nc}{\mathbb{L}}

\newcommand{\expo}{{\rm{E}}}

\newcommand{\graphs}{{\rm{Gr}}}

\newcommand{\digraphs}{{\rm{Digraph}}}

\newcommand{\pert}{{\rm{Pert}}}

\newcommand{\sh}{{\rm{Sh}}}

\newcommand{\iso}{{\rm{Iso}}}

\newcommand{\set}{{\rm{Set}}}

\newcommand{\cat}{{\rm Cat}}
\def\C{{C}}
\date{}

\begin{document}
\setlength{\baselineskip}{16pt}
\title{Super, Quantum and Non-Commutative Species}
\author{Rafael D\'\i az and Eddy Pariguan}
\maketitle
\begin{flushright}
{ Dedicated to the memory of Gian-Carlo Rota.}
\end{flushright}

\smallskip

\begin{abstract}
We introduce an approach to the categorification of rings, via the
notion of distributive categories with negative objects, and use it
to lay down categorical foundations for the study of super,
quantum and non-commutative combinatorics. Via the usual duality
between algebra and geometry, these constructions provide
categorifications for various types of affine spaces, thus
our works may be regarded as a starting point towards the
construction of a categorical geometry.\\


\end{abstract}

\section{Introduction}

This work takes part in the efforts aimed to uncover the
categorical foundations of quantum field theory QFT. The most
developed categorical approach to QFT is via the Atiyah's axioms
for topological quantum field theory \cite{MA, MA2, TU} and
closely related axioms for other theories such as conformal field theory
\cite{Se, Se2}, topological conformal field theory \cite{Cos,
kim2}, homotopy quantum field theory \cite{TU1, TU2}, and
homological quantum field theory \cite{cas3}. These categorical
approaches are most useful for the study of field theories
defined over  non-trivial topological spaces. In contrast the
study of locally defined QFT with categorical tools remains, at large, limited and conjectural
(for a promising new approach the reader may consult \cite{u}.) \\

In this work we propose an approach that, while still in its infancy, is well suited to
deal with local issues in QFT. Our approach is based on a couple
of elementary yet subtle observations. First observation: several
statements coming from QFT may be understood in a rigorous way if
regarded as taking part in formal geometry. We mention three main
examples: i) Perturbative finite dimensional Feynman integrals are
rigorous objects if one is willing to use fields which are formal
power series, i.e. power series which may be divergent,
see \cite{de} and the references therein. In some cases this
formal approach to Feynman integrals may also be applied in the
infinite dimensional context. ii) Deformation quantization \cite{bayen} of
a finite dimensional Poisson manifold $M$
consists in the construction of a star product $\star$
on the space of formal power series in a variable $\hbar$ with
coefficients in the space of smooth functions on $M$. The main
result is that the star product is essentially  determined by the
Poisson bracket on $M$. This fundamental fact was proved independently by
Fedosov \cite{fe} and De Wilde and Lecomte
\cite{dew} for symplectic Poisson manifolds;  the general case was settled by Kontsevich  \cite{kont}.
Kontsevich's work is done in the smooth context, however
when the underlying space of the Poisson manifold is
Euclidean space, his constructions applies as well in the formal
context without deep changes. Notice that in both cases, smooth
or formal,  a formal variable $\hbar$ has to be introduced. iii)
The classification of infinite dimensional Lie groups is a
notoriously difficult problem. However the classification of
simple formal supersymmetries has been accomplished by Kac
\cite{k0, k}.\\

Second observation: results in formal geometry have an underlying
categorical meaning. This basic idea was introduced in
combinatorics by Joyal  \cite{j1, j2}  via his theory of
combinatorial species, which has been developed by a
number of authors. The book \cite{Bergeron} contains a
comprehensive list of results and references in the theory of combinatorial species. Though not yet
fully appreciated by the mathematical community at large, various
constructions in the theory of combinatorial species actually has little to do
with combinatorics and may be applied as well in other categorical settings.
Considering the three formal constructions mentioned above
one is  led to the following categorical constructions: i) In Section
\ref{CFI} we  developed a categorical version of Feynman integrals.
That construction requires the extension of the notion
of combinatorial species to the notion of $G\mbox{-}C$ species,
where $G$ is a semisimple groupoid with finite morphisms and $C$ is a symmetric monoidal
category. By definition the category of $G\mbox{-}C$ species is
the category $C^G$ of functors from $G$ to $C$. ii) The
categorical version of Kontsevich star product in full generality
is given in Section \ref{kon}. The important case of constant Poisson bracket,
i.e. the categorification of the Weyl algebras, is considered in  Section \ref{qspc} where
the notion of quantum species is introduced. The
categorification of Weyl algebras allows us to look at the problem
of the normal ordering of annihilation and creation operators
\cite{DP0, DP2, RDEP3} from a new perspective. iii) The work of
Kac on the classification of formal supersymmetries opens the door
for a categorical understanding of such objects; that will be the subject of
our forthcoming work \cite{RDEP11}. Pursuing this line of research
will yield a plethora of examples of what might be called
categorical Lie algebras.\\

As the reader may guess from the previous considerations a major
requirement for this work is to have a solid
understanding of the notion of categorification.  Let us here
explain informally what do we mean by such notion, and refer
the reader to the body of this work for detailed
definitions. The notion of categorification is under
active investigation and there are various approaches to the
subject. Though implicitly present in the works of the founders of
category theory \cite{SMacLane}, the current activity on the
subject have been greatly influenced by  the works, among others,
of Baez and Dolan \cite{BaezDolan, PP}, Crane and Yetter
\cite{cy}, and Khovanov \cite{kho}. It is customary to base the foundations of
mathematics upon set theory but, as the Grothendieck's theory of topoi has shown \cite{go, SMacLane2},
in many cases it is more enlightening to look for the categorical
foundations of a given mathematical construction. The process of
uncovering the categorical foundations of a set theoretical
construction is named categorification. Let $\cat$ be the category
of essentially small categories; morphisms in $\cat(C,D)$ from a
category $C$ to a category $D$ are functors $F:C\longrightarrow
D$. Let $\set$ be the category of sets and functions as morphisms.
There is a natural functor $\cat
\longrightarrow  \set$ called  decategorification
such that:
\begin{itemize}
\item{ It sends an essentially small category $C$ to the set $\underline{C}= \ob(C)/\iso_{C}$.}
\item{ It sends a functor $F:C \longrightarrow E$
into the induced map
$$\underline{F}:\underline{C}=\ob(C)/\iso_{C}
\longrightarrow
\ob(E)/\iso_{E}=\underline{E}.$$ }
\end{itemize}
Thus  $\underline{C}$ is the set of isomorphism classes of objects in $C$. We
say that $\underline{C}$ is the decategorification of $C$ and also that $C$
is a categorification of $\underline{C}$. Notice that while a category has a
unique decategorification, a set will have many categorifications.
The motivating example, perhaps known implicitly to mankind
since its early days, is the category $\mathbb{B}$ whose objects
are finite sets and whose morphisms are bijections between finite
sets;  it is easy to check that the decategorification of
$\mathbb{B}$ is the set $\mathbb{N}$ of natural numbers.
\\

In general we are interested in the categorification of sets
provided with additional geometric or algebraic structures.  For
example one might try to find out what is the categorical analogue
of a ring. In Section \ref{jaja} we define the
categorification of a ring $R$ to be a distributive category
with negative objects provided with a $R$-valuation. The
main goal of this paper is to describe categorifications of
several types of  spaces, namely, noncommutative, quantum and super
affine spaces. This is accomplished by identifying affine spaces
with the ring of functions on them, and finding distributive categories
with natural valuations on the corresponding ring of functions.

\section{Categorification of rings}\label{jaja}

In this section we introduce the notion of categorification of
rings and provide several examples. Recall
that a monoidal category is a category $C$ provided with a
bifunctor $\odot:C\times C\longrightarrow C$ and natural
isomorphisms $\alpha_{x,y,z}: x\odot(y\odot z)\longrightarrow
(x\odot y)\odot z$ satisfying Mac Lane's pentagon
identity
$$\alpha_{x\odot
y,z,w}\alpha_{x,y,z\odot w}=(\alpha_{x,y,z}\odot 1_{w})
\alpha_{x,y\odot z,w}(1_{x}\odot\alpha_{y,z,w}).$$
A symmetric monoidal category is a monoidal category $C$ together
with natural isomorphisms $s_{x,y}: x \odot y\longrightarrow
y\odot x$ satisfying: $s_{x,y}\circ s_{y,x}=
1_x$ and Mac Lane's hexagon identity
$$(s_{x,z}\odot
1_{y})\alpha_{x,z,y}(1_{x}\odot s_{y,z})
=\alpha_{z,x,y}s_{x\odot y,z}\alpha_{x,y,z}.$$

 A categorification of a ring $R$ is a triple $(\C,N,|\ \
|)$ where $\C$ is a distributive category, $N:\C \longrightarrow \C$ is a
functor called the negative functor, and $|\ \ |: \C
\longrightarrow R$ is a map from the set of objects of $\C$ into $R$ called the valuation map. This data should satisfy the following conditions.

\begin{enumerate}
\item{ $C$ is a distributive category, i.e. $\C$ is provided with bifunctors $\oplus:\C \times \C
\longrightarrow \C$ and $\otimes:\C \times \C
\longrightarrow \C$ called sum and product, respectively.
Functors $\oplus$ and $\otimes$ are such that:
\begin{itemize}
\item{There are distinguished objects $0$ and $1$ in $C$.}
\item{The triple $(\C,\oplus,0)$ is a symmetric monoidal category with
 unit $0$.}
\item{The triple $(\C, \otimes,1)$ is a  monoidal category
with unit $1$.}
\item{Distributivity  holds. That is  for objects $x,y,z$ of $\C$ there
are natural isomorphisms $x
\otimes(y\oplus z)\simeq (x \otimes y) \oplus  (x \otimes z)$ and $(x\oplus y)\otimes z \simeq (x \otimes z) \oplus (y \otimes z)$.}
\end{itemize}

See Laplaza's works \cite{l2, l1} for a complete definition,
including coherence theorems, of a category with two monoidal
structures satisfying the distributive property.}
\item{The functor $N: \C \longrightarrow \C$ must be such that for
$x,y \in \C$ the following properties hold: $N(x\oplus
y)\simeq N(x)\oplus N(y)$,  $N(0)=0$, and $N^{2}$ is the identity functor.}
\item{The map $|\ \ |: \C \longrightarrow R$ is such that for
$x,y \in \C$ we have:
\begin{itemize}
\item{$|x|=|y|$ if $x$ and $y$ are isomorphic.}
\item{$|x\oplus y|=|x|+|y|$, $|x \otimes y| = |x||y|$, $|1|=1$, and $|0|=0$.}
\item{$|N(x)|=-|x|$.}
\end{itemize}}
\end{enumerate}

We make a few remarks regarding the notion of categorification of
rings. If $R$ is a semi-ring then a categorification of $R$ is
defined as above  omitting the existence of the functor $N.$
A categorification is said to be surjective if the valuation map is surjective.
Notice that we do not require that $\oplus$ and $\otimes$ be the
coproduct and product of $\C$, although they could be.  We stress that our definition only
demands that $|a \oplus N(a)|=0$. Demanding that $a
\oplus N(a)$ be isomorphic to $0$ would reduce drastically the scope
of our definition. In practice we prefer to write $-a$ instead of
$N(a)$. \\

A ring $R$ is  a categorification of itself, since one may consider
$R$ as the category whose object set is $R$ and with identities as the only
morphisms. The valuation map is the identity map and the negative
of $r\in R$  is $-r$. Thus, there is not existence problem
attached to the notion of categorification: all rings admit a
categorification. It will become clear from the examples given below that one should
not expect the categorification of a ring to be unique. Quite the contrary, the philosophy behind the notion of
categorification is that valuable information about
a ring  can be obtain by looking at its various categorifications, just like we can learn valuable information
about a group by looking at its various representations.\\

A functor $\varphi:C\longrightarrow D$ between distributive categories is a functor that is monoidal with respect
to $\oplus$ and $\otimes$. If both $C$ and $D$ have negative objects, then we demand in addition
that the functor   $\varphi$ respects the negative functors on $C$ and $D$. Notice that
$\otimes$ is not required to be symmetric; if it is symmetric then we say that $C$ is a symmetric distributive
category.

\begin{lem}{\em
For each distributive category $C$, there exists a
distributive category with negative objects
$\mathbb{Z}_2$-$C$, and an inclusion functor $i:C\longrightarrow
\mathbb{Z}_2\mbox{-}C$ such that for any given distributive category
with negative objects $D$ and any given functor $\varphi:C
\longrightarrow D,$ there is a unique functor $\psi:
\mathbb{Z}_2\mbox{-}C\longrightarrow D$ such
that $\psi\circ i=\varphi$, i.e. the following diagram commutes
\[
\xymatrix{ C \ar[rr]^{i} \ar[dr]_{\varphi} & &
\mathbb{Z}_2\mbox{-} C \ar[dl]^{ \psi}
\\ & D & }
\]}

\end{lem}

\begin{proof}
First define $\mathbb{Z}_2\mbox{-}C$ as the
category $\mathbb{Z}_2\mbox{-}C=C\times C$ with sums and products given by
$$(a_1,a_2)\oplus (b_1,b_2)=(a_1\oplus b_1,a_2\oplus
b_2),$$ $$(a_1,a_2)\otimes
(b_1,b_2)=(a_1\otimes b_1\oplus a_2\otimes b_2, a_1\otimes
b_2\oplus a_2\otimes b_1).$$ The negative functor $N$ is given by
$N(a,b)=(b,a)$. The inclusion functor $i:C\longrightarrow
\mathbb{Z}_2\mbox{-}C$ is given by $i(a)=(a,0).$  Given $\varphi:C
\longrightarrow D$, then the functor $\psi:
\mathbb{Z}_2\mbox{-}C\longrightarrow D$ is given by
$$\psi(a_1, a_2)=a_1 \oplus N(a_2).$$

\end{proof}

\begin{lem}\label{extn}{\em Let $|\ \ |:C\longrightarrow R$
be a valuation on a distributive category
$C$. There is a natural valuation $|
\ \ |:\mathbb{Z}_2\mbox{-}C\longrightarrow R$ on
$\mathbb{Z}_2\mbox{-}C$ such that $|i(x)|=|x|$ for  $x\in C$.}
\end{lem}
\begin{proof} Define $|\ \ |:\mathbb{Z}_2\mbox{-}C\longrightarrow R$\
\
by $|(a,b)|=|a|-|b|$.
\end{proof}

Lemma \ref{extn} allows us to define valuations with rings as codomain from
valuations with semi-rings as codomain. Next paragraphs introduce a list of
examples of distributive categories provided with valuations.\\

Let $set$ be the category of finite sets and maps as morphisms. The distributive structure on $set$ is
given by disjoint union $x\sqcup y$  and Cartesian product
$x\times y$. The map $|\ \ |:set
\longrightarrow \mathbb{N}$ sending $x$ into its cardinality $|x|$
defines a valuation on $set$.\\

Let $vect$ be the category of finite dimensional vector spaces. It
is a distributive category with  $\oplus$ and $\otimes$
defined as the direct sum and tensor product of vector
spaces. The map $|\ \ |:vect
\longrightarrow
\mathbb{N}$ given by $|V|=\dim(V)$ defines a
valuation on $vect$.\\

Let $\mathbb{Z}_{2}$-$vect$ be the category of finite dimensional
$\mathbb{Z}_2$-graded vector space. Let $V,W\in
\mathbb{Z}_{2}\mbox{-}vect$ be given by  $V=V_0\oplus V_1$ and
$W=W_0\oplus W_1$. Direct sum and tensor product on
$\mathbb{Z}_{2}\mbox{-}vect$ are given, respectively, by $V\oplus
W = (V_0\oplus W_0)\oplus (V_1\oplus W_1)$ and
$$V\otimes W=[(V_{0}\otimes W_0)\oplus (V_1\otimes W_1)]\oplus
[(V_{0}\otimes W_1)\oplus (V_1\otimes W_0)]$$ The map $|\ \
|:\mathbb{Z}_{2}\mbox{-}vect
\longrightarrow \mathbb{Z}$ given by  $|V_0\oplus
V_1|=\dim(V_0)-\dim(V_1)$ is a valuation on
$\mathbb{Z}_{2}$-$vect$.\\

Recall \cite{GCRota} that the M\"obius function $\mu:x \times x
\longrightarrow x$ of a finite partially ordered set $(x,\leq)$ is defined
as follows: for $i,k$ incomparable elements of $x$ we set
$\mu(i,k)=0.$ For $i\leq k$ the M\"obius function satisfies the
recursive relation:
$$\sum_{ i\leq j \leq k}\mu(i,j)=\left \{ \begin{array}{cc}
  1 & \mbox{if \ }i=k
 \\
  0 & \mbox{otherwise} \\
\end{array}\right.$$
Let $mposet$ be the full subcategory of the category of posets (partially ordered sets)
whose objects $(x, \leq)$ are such that each equivalence
class $c$, under the equivalence relation on $x$ generated by $\leq$,
has a minimum $m_c$ and a maximum $M_c$. The sum functor is
disjoint union of posets  and the product functor is the Cartesian
product of posets. Consider the map $|\ \ |:mposet
\longrightarrow \mathbb{Z}$ given by
$$|(x,\leq)|=\sum_{c}\mu( m_{c}, M_c),$$
where the sum runs over the set $\{ c \}$ of equivalence classes on $x$. The map $|\ \ |$ defines a
valuation on the category $mposet$.\\

Let $vman$ be the category of pairs $(M, v)$, where $M$
is a finite disjoint union of finite dimensional oriented smooth
manifolds, and $v$ is a map that sends each connected
component $c$ of $M$ into a top differential form $v(c)
\in\Omega^{\dim(c)}(c).$  Morphisms in $vman((M, v_M), (N, v_N))$ are smooth maps $f:M\longrightarrow
N$ such that $f^{\ast} v_N=v_M$. The sum functor is disjoint union and the product functor is Cartesian
product. The map  $v_{M \times N}$ sends the connected
component $c \times d$ to the differential form $$v_{M\times N}(c\times
d)=\pi_{M}^{\ast}v_M(c) \wedge\pi_{N}^{\ast}v_N(d),$$
where $\pi_{M}$ and $\pi_{N}$ are the projections from $M\times N$ onto
$M$ and $N$, respectively. The map $|\ \ |:vman
\longrightarrow
\mathbb{R}$ given by $$|(M, v)|=\sum_{c \in \pi_0(M)}\int_c
v(c)$$ defines a valuation on
$vman$ by Fubini's theorem.\\

Let $top$ the category of topological spaces with finite
dimensional $\mathbb{C}$-cohomology groups. The sum functor is disjoint
union and the product functor is Cartesian product of topological
spaces. By the K\"unneth's formula the map $|\ \ |:top
\longrightarrow \mathbb{C}[[t]]$ given by
$$|X|={\des\sum_{i=0}^{\infty}\dim_{\mathbb{C}} (H^{i}(X))t^{i}}$$
defines a valuation on $top$.\\

Let $symp$ be the category whose objects are finite dimensional
symplectic manifolds. We allow disconnected manifolds with
components of various dimension. Morphisms in $symp$ from
$(M,\omega_M)$ to $(N,\omega_N)$ are smooth maps $f: M
\longrightarrow N$ such that $f^{\ast}\omega_N=\omega_M$. The distributive structure on $symp$ is given by
disjoint union and Cartesian product, where the symplectic structure
on $M\times N$ is given by $\omega_{M\times N}=\pi_{M}^{\ast}w_M
+\pi_{N}^{\ast}w_N$. The valuation map $|\ \ |:symp \longrightarrow \mathbb{R}$ is given by
$$|(M,\omega_M)|=\sum_{c \in \pi_0(M)}\int_{c} \omega_{M}^{\frac{\dim{c}}{2}}.$$

Let $C$ be a distributive category provided with a valuation map
$|\ \ |:C \longrightarrow R$. Let $C^{\mathbb{N}}$ be the category of
$\mathbb{N}$-graded $C$-objects, i.e. the category of functors
$\mathbb{N}\longrightarrow C,$ where $\mathbb{N}$ is the
category whose objects are the natural numbers and with identities
morphisms only. The sum and product functors are given by $$(F\oplus G)(k)=F(k)\oplus
G(k),$$ $$(F\otimes G)(k)={\des \bigoplus_{i+j=k} F(i)\otimes G(j)}.$$  The  map $|\ \ |: C^{\mathbb{N}}
\longrightarrow R[[t]]$  given by $|F|={\des{\sum_{k\in
\mathbb{N}}|F(k)|t^{k}}}$ defines a valuation on
$C^{\mathbb{N}}$.\\

Let $gpd$ be the category of finite groupoids. Recall
\cite{BaezDolan} that a finite groupoid $G$ is a category such
that the objects of $G$ form a finite set, $G(x,y)$ is a finite set for
all $x,y\in G$, and all morphisms in $G$ are invertible. The sum and
product functors on $gpd$ are, respectively,  disjoint union and Cartesian product
 of categories. The valuation map $|\ \ |: gpd
\longrightarrow
\mathbb{Q}$ is given by $$|G|=\sum_{x\in \underline{G}}
\frac{1}{|G(x,x)|}.$$ This example has been exploited by D\'iaz and
Bland\'in \cite{Blan, Blan1, Blan2} in order to propose a model for the study of the
combinatorics of rational numbers.\\

Perhaps the best known example of categorification is
the following. Let $M$ be a compact topological space and $vect_M$
be the category of finite rank $\mathbb{C}$-vector bundles on $M$.
The canonical map $\pi:vect_M\longrightarrow K_{0}(M)$, where
$K_{0}(M)$ is the degree zero $K$-theory group of $M$, is a
valuation map.\\

From now on we will make the following assumption. Let $x$ be a
set of cardinality $n$, $C$ a symmetric monoidal category with product $\odot$, and $f:x\longrightarrow
C$ a map. Consider the category $\mathbb{L}(x)$  of linear
orderings on $x$. Objects in $\mathbb{L}(x)$ are bijections
$\alpha:[n]\longrightarrow x$  where $[n]=\{1,2,...,n\}$. Morphisms in $\mathbb{L}(x)$
from $\alpha$ to $\beta$ are given by:
$$\mathbb{L}(x)(\alpha,\beta)=\{\sigma:[n]\longrightarrow [n]\ |\
\beta \sigma=\alpha \
\mbox{  and}\  \sigma \ \mbox{is a bijection}\}.$$ Consider the functor $\widehat{f}:\mathbb{L}(x)\longrightarrow
C$ given by $$\widehat{f}(\alpha)=
\bigodot_{i=1}^{n}f(\alpha(i)).$$ The image $\widehat{f}(\sigma)$ of a
morphism $\sigma$ in $\mathbb{L}(x)$ is obtained using the
symmetry map of $C$. With this notation we define the
 $\odot$-product of objects in $C$ indexed by an unordered set $x$ as
follows:
$$\bigodot_{i\in x}f(i)={\rm colim}(\widehat{f}).$$
From now on we assume that our distributive
categories are such that for each map $f: x\longrightarrow C$  the colimit of
the associated functor $\widehat{f}$ exists for $\odot=\oplus$. Moreover, if
$C$ happens to be a symmetric distributive category, then we also assume
that the colimit above exist for $\odot=\otimes.$

\section{Categorification of non-commutative affine space}\label{ncesp}

In this section we begin the study of the main topic of this work, namely, the
categorification of certain affine spaces. In the previous section we
gave a precise definition of the notion of categorification of rings, and
constructed various examples. To categorify
spaces we recall the duality
$$ \mbox{\  geometry \ } \longleftrightarrow \mbox{\ algebra \ }$$
between geometry and algebra which assigns -- in its simplest
version -- to each space its corresponding ring of functions. For
example if our space is a topological space, then we consider the
ring of continuous functions on it. If instead, it is a
smooth manifold, then one considers the ring of smooth functions on it. If it is an affine
variety we consider the ring of polynomials functions on it, and so on.
This duality has been of great use in functional analysis,
algebraic geometry, non-commutative geometry, and further
applications are to be expected. The key point to keep in mind is that
once the appropriated ring of functions for a given space have
been determined, then the geometric properties of that space will
be encoded in the algebraic properties
of the corresponding ring. \\

Let $R$ be a commutative ring. The non-commutative formal
$d$-dimensional affine  space over $R$ is the space whose
associated ring of functions is $R\langle \langle x_1,\dots, x_d
\rangle
\rangle$, the ring of formal power series  with coefficients in $R$ in the non-commutative
variables $x_1,\dots,x_d$. We find a categorification of non-commutative affine
$d$-space as follows: we are going to define a distributive category $\nc_d$
such that any symmetric distributive category $C$ provided with a
valuation map $|\
\ |: C\longrightarrow R$ induces a valuation map $$|\ \
|:C^{\nc_d}
\longrightarrow  R\langle \langle x_1,\dots, x_d \rangle \rangle
.$$
\begin{defn}{\em  For $d\in \mathbb{N}^{+}$
we let $\nc_d$ be the category such that:
\begin{enumerate}
\item{ Objects of $\nc_d$ are triples $(x,\leq,f)$ where
$x$ is a finite set, $\leq$ is a linear order on $x$, and $f:x
\longrightarrow [d]$ is a map.}

\item { The set of morphisms $\nc_d((x, \leq, f), (y, \leq, g))$ is given by $$ \{   \varphi: x \longrightarrow
y  \ |  \ g\circ \varphi=f,\  \ \varphi \mbox{ is a bijection and } \  \varphi(i)<
\varphi(j)\  \mbox{ for }\  i< j \ \}.$$}
\end{enumerate}}
\end{defn}

Note that $\nc_d$ is a groupoid and that there is at most one morphism between any pair of
objects in  $\nc_d$. Given essentially small categories $C$ and
$D$, we let $D^C$ be the category of functors from $C$ to $D$.
Morphisms in $D^C(F,G)$ are natural transformations  $F\longrightarrow G$.

\begin{defn} {\em Let $C$ be a symmetric distributive category. The
category $C^{\nc_d}$ of functors from $\nc_d$ to $C$ will be
called the category of non-commutative $C$-species of type $d$. We
denote by $C_+^{\nc_d}$ the full subcategory of $C^{\nc_d}$ whose
objects are functors $F \in C^{\nc_d}$ such that
$F(\emptyset)=0$.}
\end{defn}

The category $C_+^{\nc_d}$ is most useful when $C$ does not admit
arbitrary  sums $\bigoplus_{i\in I} a_i$. Recall that the ordered disjoint
union of posets is given by $(x,\leq)\sqcup (y,\leq)=(x\sqcup
y,\leq )$, where $\leq$ on $x\sqcup y$ is such that its
restriction to $x$ agrees with the order on $x$, its restriction
to  $y$ agrees with the order on $y$, and $i \leq j$ for  $i\in
x$, $j\in y$. An ordered partition into $n$-pieces of a poset
$(x,\leq)$  is a $n$-tuple  of non-empty posets
$(x_1,\leq),\dots,(x_n,\leq)$ such that:
$$(x_1,\leq)\sqcup \dots \sqcup (x_n,\leq)=(x,\leq).$$ We
let $opar(x,\leq)$ be the set of all ordered partitions of $(x,\leq)$. To simplify notation
we denote the restrictions of  $\leq$ to the various subsets of $x$ by the
same symbol $\leq$, we hope this causes no confusion.

\begin{figure}[h]
\begin{center}
\includegraphics[width=9cm]{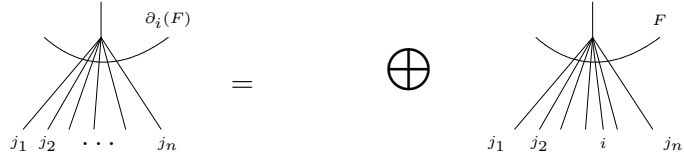}
\caption{Graphical meaning of the derivative of species. \label{fig:one}}
\end{center}
\end{figure}

\begin{defn}\label{propert}{\em Let $C$ be a symmetric
distributive category, $F,G \in C^{\nc_d},$ $G_1,\dots,G_d\in
C_+^{\nc_d}$, and $(x,\leq,f) \in \nc_d$. The following formulae
define, respectively, sum, product, composition and derivative of
non-commutative species:
\begin{enumerate}
\item{$(F + G)(x,\leq,f)=F(x,\leq,f)\oplus G(x,\leq,f)$.}
\item{$(F G)(x,\leq,f)= \bigoplus F(x_1,\leq,f{|}_{x_1})\otimes G(x_2,\leq,f{|}_{x_2})$,
where the sum runs over all pairs $(x_1,\leq)$ and $(x_2,\leq)$
such that $(x_1,\leq)\sqcup (x_2,\leq)=(x,\leq)$.}
\item{$F(G_1,\dots,G_d)(x,\leq,f)=\des{\bigoplus_{
p,g}F(p,\leq, g)\otimes \bigotimes_{b\in p} G_{p(b)}(b,\leq,
f{|}_{b})}$, where the sum runs over all $p\in opar(x,\leq)$ and
all maps $g:p\longrightarrow [d]$.}

\item{$\partial_i:C^{\nc_d}  \longrightarrow
C^{\nc_d}$ where $\partial_{i}F$ is given by
$$\partial_{i}F(x,\leq,f)={\des
\bigoplus_{\leq_{\ast}}} F(x\sqcup
\{\ast\},\leq_{\ast},f\sqcup
\{(\ast,i)\})$$
where the sum runs over all extensions $\leq_{\ast}$ of the
order on $x$ to a linear order on $x\sqcup\{\ast\}.$}
\end{enumerate}}
\end{defn}
Figure \ref{fig:one} explains the graphical meaning of the
derivative. The sum on the right hand side  runs over all possible
ways to insert the edge $(\ast,i)$ in the order set $(x,\leq,f)$.
Let $x_1,...,x_d$ be a family of non-commutative variables. Then
for $(x,\leq,f)\in \nc_d$ we write $x_{f}=x_{f(1)}\dots
x_{f(d)}$.

Let $R$ be a ring and let $C$ be a symmetric distributive category
provided with a valuation map $|\ \ |: C\longrightarrow R$.

\begin{thm} {\em  Under the conditions above $C^{\nc_d}$ is a distributive category
and the map $$|\ \ |:C^{\nc_d} \longrightarrow  R\langle
\langle x_1,\dots, x_d
\rangle \rangle $$ $$\mbox{ \ \ given for} \ \  F\in C^{\nc_d} \mbox{ \ \ by  \ \ } |F|= \des{\sum_{f:[m]\rightarrow
[d]}|F([m],f)|x_f}$$ is a valuation on $C^{\nc_d}$. Moreover we
have that
$$|\partial_i F|=\partial_i|F| \ \mbox{ and } \
|F(G_1,\dots,G_n)|=|F|(|G_1|,\dots,|G_n|).$$ }
\end{thm}
\begin{proof} For $F,G \in C^{\nc_d}$  we have that:

\

$\begin{array}{rcl} |F + G|&=&\des{\sum_{f:[m]\rightarrow
[d]}}|F\oplus G([m],f)| x_f=\des{\sum_{f:[m]\rightarrow
[d]}}|F([m],f)\oplus G([m],f)|
x_f\\

\mbox{ }&=&\des{\sum_{f:[m]\rightarrow [d]}}|F([m],f)| x_f \ +
\des{\sum_{f:[m]\rightarrow [d]}}|G([m],f)| x_f=|F|+|G|
\end{array}$

\bigskip

$\begin{array}{rcl}
|F G|&=&\des{\sum_{f:[m]\rightarrow [d]}}|F G([m],f)| x_f \\
\mbox{ }&=&\des{\sum_{f:[m]\rightarrow [d]}}\sum_{m_1\sqcup m_2=m}
|F(m_1,f{|}_{m_1})|x_{f|_{m_1}}
|G(m_2,f{|}_{m_2})| x_{f|_{m_2}}\\
\mbox{ }&=&|F| |G|.
\end{array}$

\bigskip

\noindent For $G_1,\dots,G_d \in C_+^{\nc_n}$  we have that:

\bigskip

$\begin{array}{rcl} |F (G_1,\dots,G_d)|&=&
\des{\sum_{f:[m]\rightarrow [d]}}|F
(G_1,\dots,G_d)([m],f)| x_f\\
\mbox{}&=&\des{\sum_{f:[m]\rightarrow [d]}\left| \bigoplus_{
p,g}F(p, g)\otimes \bigotimes_{b\in p} G_{p(b)}(b,
f{|}_{b})|\right|x_{f}}\\
\mbox{}&=&\des{\sum_{f:[m]\rightarrow [d]}\sum_{p,g} \prod_{b\in p}|F(p, g)|| G_{p(b)}(b,
f{|}_{b})|x_{f}}\\
\mbox{}&=&\des{\sum_{f:[m]\rightarrow [d]}\sum_{p,g} \prod_{b\in p}\left(|F(p, g)|  | G_{p(b)}(b,
f{|}_{b})|x_{f|_{b}}\right)}\\
\mbox{}&=&\des{\sum_{f:[m]\rightarrow [d]}|F([m], f)|\prod_{i=1}^{n}
\left(\sum_{g:[m_i]\rightarrow [d]}|G_i([m_i],
g_i)|x_{g_i}\right)
}\\
\mbox{}&=& |F|(|G_1|,\dots,|G_d|).\\
\end{array}$

\bigskip

$\begin{array}{rcl}
|\partial_i F|&=& \des{\sum_{f:[m]\rightarrow [d]}}|\partial_i F([m],f)| x_f\\
\mbox{}&=& \des{\sum_{f:[m]\rightarrow [d]}}\sum_{\leq_{\ast}}|  F([m]\sqcup
\{\ast\},\leq_{\ast},f\sqcup
\{(\ast,i)\})  |x_f =\partial_i|F|,\\
\end{array}$

\bigskip

\noindent since $\partial_i x_g=x_f$ if and only if the domain of $g$ is
isomorphic to $[m]\sqcup\{\ast\}$ and $g=f\sqcup\{(\ast,i)\}$.
\end{proof}

\begin{exmp}{\em Let $(x,\leq,f) \in \nc_d.$

\begin{enumerate}
\item { For $i \in [d]$ the non-commutative singleton species $X_i\in C^{\nc_d}$ is given by $X_i(x,\leq,f)=1$ if\
\ $|x|=1$ and $f(x)=i$, otherwise $X_i(x,\leq,f)=0$. We have that
$$|X_i|={\des\sum_{f:[m]\rightarrow [d]} |X_i([m],\leq,f)|x_f }=x_i .$$
where $\leq$ denotes the standard linear order on $[m]$.}

\item{The species $1\in C^{\nc_d}$ is
given by $1(x,\leq,f)=1$ if $x=\emptyset$ and $ 1(x,\leq,f)=0$
otherwise. We have that
$$|1|={\des\sum_{f:[m]\rightarrow [d]} |1([m],\leq,f)|x_f}=1.$$}
\item{The non-commutative species $NE_d\in
C^{\nc_d}$ is such that $NE_d{(x,\leq,f)}=1$ for $(x,\leq,f) \in
\nc_d.$ We have that:
$$|NE_d|={\des
\sum_{f:[m]\rightarrow[d]} |NE_d([m],\leq,f)|x_f}={\des
\sum_{f:[m]\rightarrow[d]} x_f}.$$}
\end{enumerate}}
\end{exmp}

As motivation for the study of non-commutative species we consider
the problem of finding the analogue of the notion of operads in
the non-commutative context. The reader may consult the next section for a brief summary on
operads. Notice that our definition of non-commutative operads is actually
an analogue of the notion of non-symmetric operads. For the next proposition, and in
other similar situations, we regard $[d]$ as the category with
identity morphisms only.

\begin{prop}{\em Let $C$ be a symmetric distributive category, then
$(C_+^{[d]\times \nc_d},\circ,(X_1,\dots,X_d))$ is a monoidal
category where for $(F_1,\dots,F_d) ,(G_1,\dots,G_d)\in
C_+^{[d]\times \nc_d}$ and $i\in[d]$ we set:
$$(F_1,\dots,F_d)\circ
(G_1,\dots,G_d)=(F_1(G_1,\dots,G_d),..., F_i(G_1,\dots,G_d),..., F_n(G_1,\dots,G_d)).$$}
\end{prop}

Let $(C,\odot,1)$ be a monoidal category. A monoid $M$ in $C$
\cite{SMacLane, SMacLane2} is an object $M\in C$ together
with morphisms $m\in C(M\odot M, M)$ and $u\in C(1,M)$ such that
the following diagrams commute

\[
\begin{array}{ccc}

\xymatrix{ M\odot M \odot M \ar[r]^{\mbox{  } \  \  1\odot m} \ar[d]_{m\odot 1 } & M\odot
M
\ar[d]^{m}
\\ M\odot M  \ar[r]^{m} & M }
 &\hspace{1cm} &\xymatrix{ 1 \odot M \ar[r]^{u\odot 1} \ar[dr] &
M\odot M \ar[d]^{m} & M\odot 1 \ar[l]_{1\odot u}
\ar[dl]
\\ \mbox{} & M  & \mbox{}}
\end{array}
\]
{\noindent}where the diagonal arrows are the canonical isomorphisms coming from
the properties of the unit element in a symmetric monoidal category.

\begin{defn}{\em
A non-commutative $d$-operad in $C$ is a monoid  in
$(C_+^{[d]\times \nc_d},\circ,(X_1,\dots,X_d))$.}
\end{defn}

Our next result gives an explicit description of non-commutative
operads.

\begin{thm}{\em
A non-commutative $d$-operad in $C$ is given by a collection of
objects $O=\{O^{j}_f\}$ where  $O^{j}_f \in C,$ $j\in[d],$ and
$f:[m]\longrightarrow [d]$ is a map with domain $[m]$ for some
$m\geq 1$. In addition there should be unit maps
$u_j:1\longrightarrow O^{j}_{1_j}$ and composition maps
$$\gamma: O^{j}_f \otimes \bigotimes_{i=1}^{m}
O^{f(i)}_{g_i} \longrightarrow O^{j}_{\sqcup g_i},$$ where $
g_i:[k_i]\longrightarrow [d]$, and $\sqcup g_i:
\sqcup [k_i]\longrightarrow [d]$ is  given by
$${\des (\bigsqcup g_i)|_{[k_j]}=g_j.}$$  These data should satisfy the associativity axiom:

{\small
\[\xymatrix @C=.2in { {\des O^{j}_f\otimes
\left(\bigotimes_{i}O^{f(i)}_{g_i}\right)\otimes \left(\bigotimes_{s}
O^{\sqcup_{i}g_i(s)}_{h_{is}}\right)}
 \ar[rr]^{\ \mbox{ }\ \ \ \gamma \otimes \id}
\ar[dd]_{\mbox{shuffle}} & & {\des O^{j}_{\bigsqcup_{i} g_i}\otimes
\bigotimes_{s=1}^{|k|} O^{\sqcup_{i}g_i(s)}_{h_{is}}}
\ar[d]^{\gamma} \\
& & {\des O^{j}_{\bigsqcup_{i,s}h_{is}}}  \\
{\des O^{j}_f\otimes \bigotimes_i \left(
O^{f(i)}_{g_i}\otimes\bigotimes_{s=1}^{k_i}
O^{g_i(s)}_{h_{is}}\right)}
 \ar[rr]^{\mbox{ }\ \ \ \ \id \otimes (\otimes_{m}\gamma )} & &
{\des O^{j}_f\otimes \left(\bigotimes_{i=1}^{m}
O^{f(i)}_{\sqcup_{s}h_{is}}\right)} \ar[u]_{\gamma} }\]}

\noindent Let $1_j:[1]\longrightarrow [d]$ be given by $1_j(1)=j$ for
$j\in[d]$. The following unity axioms must hold:

\[
\xymatrix{  O^{j}_f = O^{j}_f
\otimes 1\otimes\dots\otimes 1 \ar[r]^{\ \ \ \  1} \ar[d]_{1 \otimes u_{f(1)} \otimes \dots \otimes u_{f(1)}}  &
{\des O^{j}_f= O^{j}_{ \sqcup 1_{f(i)}}}
\\ {\des O^{j}_f
\otimes \bigotimes_{i=1}^{m} O^{f(i)}_{1_{f(i)}}}\ar[ur]_{\gamma}  }
\]

\[
\xymatrix{ 1\otimes O^{j}_f = O^{j}_f
 \ar[r]^{\ \ \ \ \ \  1} \ar[d]_{u_j}  &
\des O^{j}_f
\\ O^{j}_{1_j}\otimes
O^{j}_f\ar[ur]_{\gamma} }
\]}
\end{thm}

\begin{figure}[!h]
\begin{center}
\includegraphics[width=10cm]{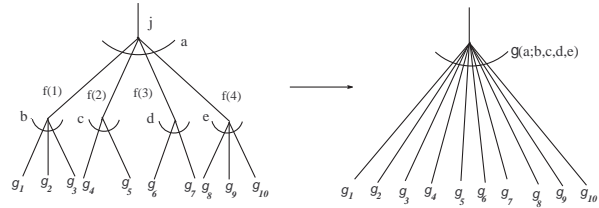}
\caption{Example of an application of the  map $\gamma$. \label{fig:multipli}}
\end{center}
\end{figure}

Figure
\ref{fig:multipli} illustrates the meaning of the composition maps
$\gamma$.

\begin{figure}[!h]
\begin{center}
\includegraphics[width=6cm]{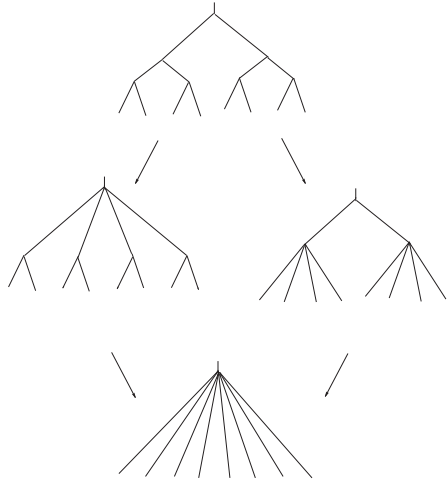}
\caption{Example of  the associativity for $\gamma$. \label{fig:asocia}}
\end{center}
\end{figure}

Figure \ref{fig:asocia} explains graphically the associativity
axiom. Next example provides a simple construction, a non-commutative analogue of the endomorphisms operad, that shows that there are plenty of non-commutative operads.

\begin{exmp}{\em
Let $C$ be a monoidal category. There is a non-commutative
$d$-operad $E$ associated with any $d$-sequence $(a_1,\dots,a_d)$
of objects  in $C$ given for $i\in [d]$ and $f: [m]\longrightarrow
[d]$ by
$$E^{i}_f=C(\bigotimes_{j \in [m]} a_{f(j)},
a_j).$$ }
\end{exmp}

\begin{exmp}{\em If $F$ is a non-commutative species, then we let $F_+$ be the species
such that $F_+((x,\leq,f))=F((x,\leq,f))$ if $x$ is non-empty and $F_+(\emptyset)=0$.
 The $n$-tuple $(NE_{d,+},...,NE_{d,+})$ is a non-commutative
$n$-operad. }
\end{exmp}

Next we shall define a non-commutative analogue for the binomial
coefficients \cite{GCRota}. Let $R$ be a commutative ring.

\begin{defn}{\em \begin{enumerate}
\item{ Consider a family $\{s_{f}\},$ where $s_{f}:\mathbb{N}\longrightarrow
R$ is a map, and $f$ is a map from some $[m]$ into $[d]$. We called such a family a non-commutative
multiplicative sequence if it is such that for $a,b\in\mathbb{N}$
the following identity holds: $$s_f(a+b)=
\sum_{i=1}^{m}s_{f_{< i}}(a)s_{f_{\geq i}}(b),$$
where $f_{< i}:[1,i-1]\longrightarrow [d]$ and $f_{\geq
i}:[i,m]\longrightarrow [d]$ are the restrictions of $f$ to the
appropriated domains. }

\item{Consider a family
$\{s_{f,i}\}$, where $s_{f,i}:\mathbb{N}\longrightarrow
R$ is a map, $f:[m]\longrightarrow [d]$ is a map with domain $[m]$, for some $m \in \mathbb{N}$,
and $i \in [d].$ We call such a family a non-commutative
compositional sequence if for $a,b\in\mathbb{N}$ the following
identity holds:
$$s_{f,i}(a+b)={\des \sum_{p} s_{p,i}(a)\prod_{j=1}^{k}
s_{f_j,p(j)}(b)},$$ where the sum runs  the ordered partitions
$$(x_1, \leq, f_1)\sqcup\dots \sqcup(x_k, \leq, f_k)=([m], \leq, f) \mbox{\ and \ the \ maps \ }
p:[k]\longrightarrow [d].$$}
\end{enumerate}}
\end{defn}

The following result gives a simple construction that generates
non-commutative multiplicative sequences, and provides a
categorical interpretation for it. Let $s={\des
\sum_{f:[m]\rightarrow [d]} s_f x_f}$ be a non-commutative formal power series. For  $a\in \mathbb{N}$ we
set $s^{0}=1,$ $s^{a+1}=s^{a}s,$ and  $$s^{a}={\des
\sum_{f:[m]\rightarrow [d]} s_f (a)x_f}.$$ For $S\in C^{\nc_d}$ we define recursively $S^0=1$,
 and $S^{a+1}=S^{a}S$.

\begin{thm}{\em
\begin{enumerate}
\item{  For $s\in R\langle\langle
x_1,\dots,x_d\rangle\rangle$ the family $\{s_f\}$ defined above is
a noncommutative multiplicative sequence. }
\item{ Let $S\in C^{\mathbb{L}_d}$ be such that $|S|=s$, for
$a\in\mathbb{N}^{+}$ let $opar([m],a)$ be the set of ordered
partitions of $[m]$ into $a$ blocks. Then we have that
$$s_{f}(a)=\left| \bigoplus_{\pi \in opar([m],a)}
\bigotimes_{b \in \pi}S(b,\leq,f|_b)\right|.$$}
\item{ Let $S\in C^{\mathbb{L}_d}$ be given by $S=1 - F$ where $F \in C^{\mathbb{L}_d}_+.$ The non-commutative species $S^{-1}$ sending the empty set into $1$ and a non-empty linearly order colored set $(x, \leq, f)$ into
    $$S^{-1}(x, \leq, f)= \bigoplus_{\pi \in opar(x, \leq)}
\bigotimes_{b \in \pi}F(b,\leq,f|_b)$$
is such that $$|S^{-1}||S|=1=|S||S^{-1}| .$$    }
\end{enumerate}}
\end{thm}
\begin{proof}
We proof the second part.
$$ \sum_{f:[m] \rightarrow [d]} s_{f}(a) x_f= |S|^{a}=|S^{a}|
= \sum_{f:[m] \rightarrow [d]} |S^{a}([m],\leq,f)| x_f=$$ $$=
\sum_{f:[m] \rightarrow [d]} \left| \bigoplus_{\pi\in opar([m],a)}
\bigotimes_{b \in opar}S(b,\leq,f|_b)\right|
x_f.$$
\end{proof}

Notice that we are not claiming that $S S^{<-1>}=1.$ It would be nice to have such an identity, but
in most cases we have to deal with the  weaker identity  $|S||S^{<-1>}|=1.$\\

\begin{figure}[h!]
\begin{center}
\includegraphics[width=2.5in]{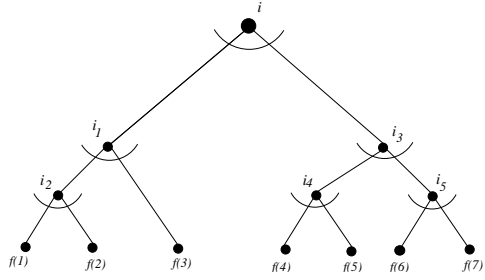}
\caption{ Example of a colored planar tree.\label{fig:tree}}
\end{center}
\end{figure}

We let $\digraphs$ be the category whose objects are directed
graphs. A direct graph is a triple $(V,E,(s,t))$ where
$V$ is the set of vertices, $E$  is the set of edges and $(s,t):E
\longrightarrow V\times V$ is a map.  A rooted tree
is a direct graph  such that there is a distinguished
vertex $r$ called the root; and for each vertex $v$ there is a unique
directed path from $v$ to $r$. Notice that necessarily
$|s^{-1}(r)|=0$. A vertex $v$ is a leaf if $|t^{-1}(v)|=0$. A
planar rooted tree is a tree together with a linear order on
$t^{-1}(v)$ for each vertex $v$. We shall consider colored or labeled planar
trees, i.e. planar trees provided with a map $l:V\longrightarrow
[n]$. Figure \ref{fig:tree} shows an example of a colored planar
rooted tree. Let $\mathrm{T}$ be the category whose objects
are colored planar rooted trees, and whose morphisms are morphisms
between the underlying directed graphs that preserve labels
and the linear orders associated with each vertex. A colored set is a set $x$ together with
a map $f:x \longrightarrow [n]$. Given a
linear ordered colored set $(x, \leq, f)$ and  $a \in\mathbb{N}$, we let
$\mathrm{T}^{a}_i(x)$ be the full subcategory of $\mathrm{T}$
whose objects are colored planar rooted trees $\gamma$ such that:
\begin{itemize}
\item $l(r)=i$.

\item The set of leaves is $x$ and  $l(i)=f(i)$ for $i\in x$.

\item The linear order on $x$ agrees with the order induced by
the planar structure on $\gamma$.
\end{itemize}

Also we define $\mathrm{T}^{a}_i(x)$ as the full subcategory of
$\mathrm{T}_i(x)$ whose objects are such that any path from a leave
to the root has length $a$. \\

Next result provides a source of non-commutative compositional
sequences, and also  a categorical interpretation for it. Moreover
we show that, in a sense made clear in the statement of the next
theorem, most $C$-valued non-commutative species have a compositional inverse.
Let $s=(s_1,\dots,s_d)\in R\langle\langle
x_1,\dots,x_d\rangle\rangle^{d}$ be such that $s(0)=0$ and
$\partial_js_i=\delta_{ij}.$ For $a \in \mathbb{N}$ we  set
$s^{<0>}=x=(x_1,...,x_d),$ and $s^{<a+1>}=s^{a} \circ s.$  We also
set $$s^{<a>}_i={\des
\sum_{f:[m]\rightarrow [d]} s_{f,i}(a)  x_f}.$$
For $S\in C_+^{[d]\times \nc_d}$  we define recursively
$S^{\langle 0\rangle}=(X_1,...,X_d)$ and $S^{\langle
a+1\rangle}=S^{\langle a\rangle}\circ S$.

\begin{thm}{\em
\begin{enumerate}
\item{The sequence $\{ s_{f,i}\}$ defined above
is a non-commutative compositional sequence. }
\item{Suppose that $S\in C_+^{[d]\times \nc_d}$ satisfies $|S|=s$, then
$${\des s_{f,i}(a)=\left|\bigoplus_{\gamma\in \underline{T_i^{a}(x)}}
\bigotimes_{v\in V_{\gamma}\setminus x} C_{l(v)}(t^{-1}(v), \leq, l)\right|,}$$ where $l$ denotes the coloring
of the graph $\gamma$.}
\item{If $S\in C_+^{[d]\times \nc_d}$ is such that $S_i=X_i-F_i$, where
$F_i \in C_+^{[d]\times \nc_d}$ and  $F_{i}(x,\leq,f)=0$ for
$|x|\leq 1$, then the non-commutative species $S^{<-1>} \in
C_+^{[d]\times\nc_d}$  given by
$${\des S^{<-1>}_{i}(x, \leq,f)=\bigoplus_{\gamma\in \underline{T_i(x)}}
\bigotimes_{v\in V_{\gamma}\backslash x} F_{l(v)}(t^{-1}(v), \leq, l)},$$
is such that
$$|S|\circ |S^{<-1>}|=(X_1,...,X_d)= |S^{<-1>}| \circ |S|.$$}
\end{enumerate}}
\end{thm}

\section{Categorification of affine space}\label{cmmsp}

We begin this section recalling the construction of the category
of $C$-valued commutative species, following the approach
introduced by Joyal \cite{j1, j2}, and a fully developed by
Bergeron, Labelle and Leroux \cite{Bergeron}. The notion of
species, under the name of collections, has also appeared in
algebraic topology, for example, in the works of Boardman and Vogt
\cite{bv}. We show that commutative and non-commutative
$C$-species are intertwined by a pair of adjoint functors. Let
$\mathbb{B}^{d}$ be the category whose objects are pairs $(x,f)$
where $x$ is a finite set and $ f: x
\longrightarrow [d]$ is a map. Morphisms in $\mathbb{B}^{d}$ are given by $$\mathbb{B}^{d}((x,f),(y,g))=\{ \alpha:
x\longrightarrow y \ | \
\alpha \ \mbox{is bijective and} \ \ g\circ \alpha=f\}.$$

\begin{defn}{\em
The category $C^{\mathbb{B}^{d}}$ of  functors from
$\mathbb{B}^{d}$ to $C$ is called the category of $C$-species of
type $d$. We let $C_+^{\mathbb{B}^{d}}$ be the full subcategory of
$C^{\mathbb{B}^{d}}$ whose objects are functors $F$ such that
$F(\emptyset)=0$.}
\end{defn}

Consider the species $par:\mathbb{B} \longrightarrow set$ that sends a finite
set $x$ into the set of all its partitions, i.e. families of non-empty disjoints subsets of
$x$ with union equal to $x$.

\begin{defn}\label{propi}{\em
Let $C$ be a symmetric distributive category. Let $F,G\in
C^{\mathbb{B}^{d}}$ and $G_1,\dots,G_d\in C_+^{\mathbb{B}^{d}}$
and $(x,f)\in \mathbb{B}^{d}$, the following formulae defines sum,
product, composition and derivative for commutative species:

\begin{enumerate}
\item{$(F + G)(x,f)=F(x,f)\oplus G(x,f)$.}
\item{$FG(x,f)={\des\bigoplus_{y\subset x}F(y,f{|}_{y})}\otimes
G(x-y,f{|}_{x-y}).$}
\item{$F(G_1,\dots,G_d)(x,f)=\des{\bigoplus_{p,g}F(p,g)\otimes
\bigotimes_{b\in p}G_{p(b)}(b,f{|}_{b}) }$, where
$p\in par(x)$ and $g:p\longrightarrow [d]$ is a map.} \item{
$\partial_i: C^{\mathbb{B}^{d}}  \longrightarrow
C^{\mathbb{B}^{d}}$ where $\partial_{i}F$ is given by the formula
$\partial_i F(x,f)=F(x
\sqcup\{\ast\}, f\sqcup
\{(\ast,i)\}).$ }
\end{enumerate}}
\end{defn}

For $(a_1,\dots,a_d)\in\mathbb{N}^{d}$ we write $a!=a_1!
\dots a_d!$, $[a]=([a_1],\dots,[a_d])$ and $x^{a}=x_1^{a_1}\dots x_d^{a_d}$.
Let $R$ be a ring of characteristic $0$ and $C$ be a symmetric monoidal category
provided with a valuation map $|\ \ |:C \longrightarrow R$.

\begin{thm}{\em  The
map $|\ \ |: C^{\mathbb{B}^{d}}  \longrightarrow R[[x_1,x_2,\dots,
x_d]]$ given by $$|F|={\displaystyle
\sum_{a \in \mathbb{N}^d}|F[a]|\frac{x^{a}}{a!}}$$
is a valuation on $C^{\mathbb{B}^{d}}$. Moreover $|\partial_i
F|=\partial_i|F|$ and $|F(G_1,\dots,G_d)|=|F|(|G_1|,\dots,|G_d|),$
for $F\in C^{\mathbb{B}^{d}}$ and $G_1,\dots,G_d \in
C^{\mathbb{B}^{d}}_+.$ }
\end{thm}

\begin{exmp} {\em Let $(x,f)\in \mathbb{B}^{d}$.
\begin{itemize}
\item{ The singleton specie $X_i\in
C^{\mathbb{B}^{d}},$ for $i\in [d],$ is such that $X_i(x,f)=1$ if
$|x|=1$ and $f(x)=i$, otherwise $X_i(x,f)=0$. }
\item{The species $1\in C^{\mathbb{B}^{d}}$ is given by $1(x,f)=1$ if $x=\emptyset$, and $1(x,f)=0$ otherwise.
} \item{The exponential species $E \in C^{\mathbb{B}^{d}}$ is
given by $E{(x,f)}=1.$}
\end{itemize}

\noindent It should be clear that $|X_i|=x_i$, $|1|=1$ and
$|\expo|=e^{x_1+\dots+x_n}.$}
\end{exmp}

\begin{exmp}{\em Let us give an example of a combinatorial species with a biological flavor. Let
$ADN \in set^{\mathbb{B}^{4}}$ be the species such that for
$A,T,C,G\in \mathbb{B}$ we have:
\begin{enumerate}

\item{$ADN(A,T,C,G)$ is the set of
ordered restricted matchings on $A\sqcup T\sqcup C\sqcup G$.}
\item{An ordered restricted matching $\alpha$ is a map
$\alpha:\{0,1\}\times [n]\longrightarrow A\sqcup T\sqcup C\sqcup
G$ such that:
\begin{itemize}
\item{$\alpha$ is a bijection.}
\item{$\alpha(0,i)\in A\ ( \mbox{respectively } C)$ if and only if $\alpha(1,i)\in T\  ( \mbox{respectively } G)$.}
\item{$\alpha(1,i)\in A\ ( \mbox{respectively } C)$ if and only if $\alpha(0,i)\in T\  ( \mbox{respectively } G)$.}
\end{itemize}
\item{It  is not hard to check that the valuation $|\adn|\in\mathbb{N}[[a,t,c,g]]$  of $\adn$ is given by
$$|\adn|(a,t,c,g)={\des \frac{1}{1-2at-2cg}}.$$}}

\end{enumerate}}

\end{exmp}

Recall  \cite{KM, May3} that an operad $O$ in a symmetric monoidal
category $(C,\odot ,1)$ consists of a family $O=\{O_d
\}$, with $d \in \mathbb{N}$ and $O_d \in C,$ together with the following structural maps:
\begin{enumerate}
\item{A composition law $\gamma : O_k \odot O_{j_{1}} \odot \cdots \odot
O_{j_{k}}\longrightarrow O_j$, for $k\geq 1$ and $j_{s} \geq 0$
such that $\sum_{s=1}^k j_{s}=j$.}
\item{A right action $O_d \times S_{d} \longrightarrow O_d$ of the symmetric group $S_d$ on $O_d.$}
\item{An unit map $\eta:1 \longrightarrow O_1$.}
\end{enumerate}

The structural maps are required to be associative, unital and
equivariant in the appropriated sense \cite{gk, KM}. Non-symmetric operads are defined omitting
the actions of the symmetric groups. Let $\cat$
be the category of essentially small categories.

\begin{prop}{\em
The collection $\{C_+^{\mathbb{B}^{d}}\}_{d\geq 0}$ is an operad
in $\cat$.}
\end{prop}
\begin{proof}
For  $k\geq 1$ the composition map
$$\gamma_k:C^{\mathbb{B}^{k}}\times
C^{\mathbb{B}^{n_1}}\times\dots\times
C^{\mathbb{B}^{n_k}}\longrightarrow
C^{\mathbb{B}^{n_1+\dots+n_k}}$$ is given by
$$\gamma_k(F,G_1,\dots,G_k)=F(G_1,\dots,G_k).$$ Given $\sigma\in
S_d$ and $F\in C^{\mathbb{B}^{d}},$ let $F\sigma\in
C^{\mathbb{B}^{d}}$ be such that $F\sigma(x,f)=F(x,\sigma f)$. The
map $$C^{\mathbb{B}^{d}}\times S_d\longrightarrow
C^{\mathbb{B}^{d}}$$ given by $(F, \sigma) \longmapsto F\sigma$
provides $C^{\mathbb{B}^{d}}$ with a $S_d$-action.

\end{proof}

For the next result we assume the sum functor $\oplus$ on $C$
behaves as a coproduct, i.e. any morphisms
$\varphi:\bigoplus_i c_i \longrightarrow d$  in $C$ is uniquely determined
by a family of morphisms $\varphi:c_i \longrightarrow d.$

\begin{thm}\label{adjf}{\em  Let $C$ be a symmetric distributive category.
\begin{enumerate}
\item{
The following maps are functorial:
\begin{itemize}
\item{ $S:C^{\mathbb{B}^{d}}\longrightarrow
C^{\nc_d}$ given by  $SG(x,\leq,f)=G(x,f)$.} \item{ $\Pi:C^{\nc_d}
\longrightarrow C^{\mathbb{B}^{d}}$ given by $\Pi
F(x,f)={\des\bigoplus_{\leq}} F(x,\leq,f)$ where the sum runs over
the linear orders $\leq$ on $x$.}
\end{itemize}}
\item{$\Pi$ is a left adjoint of $S$.}
\end{enumerate}}
\end{thm}

\begin{proof}
Let $F \in C^{\nc_d}$ and $G \in C^{\mathbb{B}^{d}}$. We must show
that $C^{\mathbb{B}^{d}}(\Pi F,G)=C^{\nc_d}(F,SG)$. An element
$S\in C^{\mathbb{B}^{d}}(\Pi F,G)$ is a natural transformation $
S:\Pi F\longrightarrow G$. $S$ is given by  morphisms $S(x,f):
\Pi F(x,f)\longrightarrow G(x,f),$ one for each pair $(x,f)\in
\mathbb{B}^{d},$ such that for any morphisms $\varphi:(x,f)
\longrightarrow (y,g)$ in $\mathbb{B}^{d}$ the diagram
\[
\xymatrix{ \Pi F(x,f) \ar[r]^{S(x,f)} \ar[d]_{\Pi F(\varphi)} & G(x,f)
\ar[d]^{G(\varphi)}
\\ \Pi F(y,g) \ar[r]^{} & G(y,g) }
\]
is commutative. A morphism ${\des S(x,f):\bigoplus_{\leq}F(x,\leq,f)
\longrightarrow G(x,f)}$ is the same as a  family of morphisms
$S(x,\leq,f):F(x,\leq,f)\longrightarrow SG(x,\leq,f),$ which in turns
defines a natural transformation $S:F\longrightarrow SG$.
\end{proof}

Our next result gives and interpretation in terms of formal power
series for the couple of adjoint functors given above.  Consider
the $R$-linear  map $\pi:R\langle\langle
x_1,\dots,x_d\rangle\rangle\longrightarrow R[[x_1,\dots,x_d]]$
given on monomials by $$\pi(x_{f})={\des
\prod_{i=1}^{n} x_i^{|f^{-1}(i)|}}$$ for $x_{f}\in R\langle\langle
x_1,\dots,x_d\rangle\rangle$. Consider also the linear map
$s:R[[x_1,\dots,x_d]]\longrightarrow R\langle\langle
x_1,\dots,x_d\rangle\rangle$ given by $${\des
s\left(\frac{x^{a}}{a!}\right)=\sum_{f}\prod_{i\in
[a_1+\dots+a_d]} x_{f(i)}},$$ where the sum runs over the maps
$f:[a_1+\dots+a_d]\longrightarrow [d]$ such that $f^{-1}(i)=a_i$.

\begin{thm} {\em Let $C$ be a symmetric distributive category provided with
a $R$-valuation. The following diagrams commute:
\[
\begin{array}{ccc}
\xymatrix{C^{\mathbb{B}^{n}} \ar[r]^{|\mbox{ }|\ \ \ \ \ \ \   } \ar[d]_{S} & R[[x_1,\dots,x_d]]
\ar[d]^{s}
\\ C^{\nc_n} \ar[r]^{|\mbox{ }|\ \ \ \ \ \ \ \  } & R\langle\langle
x_1,\dots,x_d\rangle\rangle } & \hspace{2cm}&
\xymatrix{ C^{\nc_n} \ar[r]^{|\mbox{ }| \ \ \ \ \ \ \ } \ar[d]_{\Pi } & R\langle\langle
x_1,\dots,x_d\rangle\rangle
\ar[d]^{\pi}
\\ C^{\mathbb{B}^{d}} \ar[r]^{|\mbox{ }|\ \ \ \ \ \ \ } & R[[x_1,\dots,x_d]] }

\end{array}
\]}

\end{thm}
\begin{proof} For $F\in C^{\mathbb{B}^{d}}$ we have that
\begin{eqnarray*}
|SF|&=&{\des\sum_{f:[m]\rightarrow [d]}|SF([m],\leq,f)|x_f}=
{\des\sum_{f:[m]\rightarrow [d]}|F([m],f)|x_f}\\
\mbox{} &=&{\des\sum_{a\in \mathbb{N}^{d}}\sum_{f:[m]\rightarrow
[d],\  f^{-1}(i)=a_i}|S(F)([m],f)|x_f}\\
\mbox{}&=& \des{\sum_{a\in
\mathbb{N}^{d}}F[a]\sum_{f:[m]\rightarrow
[d],\  f^{-1}(i)=a_i}\prod_{i\in [a_1+\dots+a_d]} x_{f(i)}}={\des
\sum_{a\in
\mathbb{N}^{d}}|F[a]|s(\frac{x^{a}}{a!})=s|F|.
}
\end{eqnarray*}

\noindent Let $F\in C^{\nc_d}$ then we have that:
\begin{eqnarray*}
\pi|F|&=&{\des\sum_{f:[m]\rightarrow [d]}F([m],\leq,f)\pi(x_f)=
\sum_{a\in \mathbb{N}^{d}}\left( \sum_{f^{-1}(i)=a_i} |F([m],\leq,f)|\right) \prod_{i=1}^{d} x_i^{|f^{-1}(i)|}}\\
\mbox{}&=&{\des \sum_{a\in \mathbb{N}^{d}}\left(
\sum_{|f^{-1}(i)|=a_i}|F([m],\leq,f)| \right)a!\frac{x^{a}}{a!}}\\
\mbox{}&=&{\des \sum_{a\in
\mathbb{N}^{d}}\sum_{\sigma \in S_{[a_1+\dots+a_d]}}| F( [a_1+\dots+a_d],\leq,f\circ
\sigma)}\\
\mbox{}&=&{\des \sum_{a\in
\mathbb{N}^{d}}\left( \sum_{\leq}|F([a_1,\dots,a_d],\leq,f)|\right)\frac{x^{a}}{a!}}=|\Pi F|.
\end{eqnarray*}

\end{proof}

Next we define the shuffle bifunctor which intertwines the
product on commutative and non-commutative species.

\begin{defn}{\em
The shuffle bifunctor $\sh: C^{\nc_d}\times
C^{\nc_d}\longrightarrow C^{\nc_d}$ is given on objets by
$$\sh(F,G)(x,\leq,f)={\des \bigoplus_{y\subset x}
F(y,\leq,f|_{y})\otimes G(x-y,\leq,f|_{x-y})}$$ and is defined on morphisms in the natural way.}
\end{defn}
\begin{thm}{\em \begin{enumerate}\item{ For  $F,G\in C^{\nc_d}$
the functor $\Pi$ satisfies:
\begin{enumerate}
\item{$\Pi(F + G)=\Pi (F)+ \Pi(G)$.}
\item{$\Pi(FG)=\Pi(F)\Pi(G)$.}
\item{$\Pi\partial_i F=\partial_i \Pi F$.}
\end{enumerate}}
\item{For $F,G\in C^{\mathbb{B}^{d}}$
the functor $S$ satisfies:
\begin{enumerate} \item{$S(F + G)=S(F) + S(G)$ for $F,G\in
C^{\mathbb{B}^{d}}$.}
\item{$S(FG)=\sh(SF,SG)$.}
\item{$\partial_{i}SF=S\partial_i
F$.}
\item{$\Pi S=d!$.}
\end{enumerate}}
\end{enumerate}}

\end{thm}
\begin{proof}Let $F,G\in C^{\nc_d}$ then we have that:
\[
\begin{array}{lcl}
\Pi(F + G)(x,f)&=&{\displaystyle\bigoplus_{\leq} F\oplus
G(x,\leq,f)=\bigoplus_{\leq} F(x,\leq,f)\oplus G(x,\leq,f)}\\
\mbox{}&=&(\Pi F\oplus
\Pi G)(x,f).\\
\Pi(FG)(x,f)&=&{\displaystyle \bigoplus_{\leq} FG(x,\leq,f)=
\bigoplus_{\leq}\bigoplus_{x_1\sqcup x_2=x} F(x_1,\leq,f|_{x_1})\otimes G(x_2,\leq,f|_{x_2})}\\
\mbox{} &=&\Pi F\Pi G(x,f).\\
\Pi\partial_i F(x,f)&=& {\displaystyle
\Pi F(x \sqcup\{\ast\}, f\sqcup
\{(\ast,i)\})=\bigoplus_{\leq_{\ast}}F(x \sqcup\{\ast\},\leq_{\ast}, f\sqcup
\{(\ast,i)\}) }\\
\mbox{}&=& \partial_{i}\Pi F(x,f).
\end{array}
\]
Let $F,G\in C^{\mathbb{B}^{d}}$ then we have that:
\[
\begin{array}{lcl}
S(F + G)(x,\leq,f)&=& F\oplus G(x,f)= F(x,f)\oplus G(x,f)=(SF\oplus
SG)(x,\leq,f).\\
\ \ \ \ S(FG)(x,\leq,f)&=& FG(x,f)=
{\displaystyle\bigoplus_{y\subset x} F(y,f|_{y})\otimes G(x\setminus y,f|_{x\setminus y})}\\
\mbox{} &=&{\displaystyle\bigoplus_{y\subset x} SF(y,\leq, f|_{y})\otimes
SG(x\setminus y,\leq, f|_{x\setminus y})}\\
\mbox{} &=& \sh(SF,SG)(x,\leq,f).\\
\ \ \ \ \ S\partial_i F(x,\leq,f)&=& {\displaystyle \bigoplus_{\leq_{\ast}}
 SF(x \sqcup\{\ast\},\leq_{\ast}, f\sqcup
\{(\ast,i)\})=F(x \sqcup\{\ast\}, f\sqcup
\{(\ast,i)\}) }\\
\mbox{}&=& S\partial_{i}F(x,\leq,f).\\
 \ \ \ \ \ \ \ \ \Pi SF (x,f) &=& \bigoplus_{\leq}SF(x,\leq,f)=\bigoplus_{\leq} F(x,f)=d! F(x,f).
\end{array}
\]

\end{proof}

\begin{defn}{\em
\begin{enumerate}
\item{ Consider a sequence  $\{s_{n}\}_{n\in\mathbb{N}^{d}}$ where $s_{n}:\mathbb{N}\longrightarrow
R$ is a map.  We call such a sequence a  polybinomial sequence if for
$a,b\in\mathbb{N}$ we have: $$s_n(a+b)={\des \sum_{j\in
\mathbb{N}^{d}} {n_1 \choose j_1} \dots
{n_d \choose j_d} s_j(a)s_{n-j}(b)}.$$ }
\item{Consider a sequence of maps
$\{s_{n,i}\}_{n\in\mathbb{N}^{d}}$, where
$s_{n,i}:\mathbb{N}\longrightarrow R$ is a map and  $i\in [d]$. We call such a sequence a  polymultinomial sequence if it
is such that
for $a,b\in\mathbb{N}$ the following identity holds:
\begin{equation*}
s_{n,i}(a+b)={\des \sum_{m,p}\frac{s_{m,i}(a)}{m!}\prod_{l=1}^d {n_l \choose \overline{p}_l } \prod_{j=1}^d
\prod_{k=1}^{m_j}s_{\overline{p}_{j,k},j}(b)   }
\end{equation*}
where $m=(m_1,...,m_d) \in \mathbb{N}^d,$ and $p$ is a map that sends triples $l,j,k$ such that
$l,j \in [d]$ and $k \in [m_j]$ into $\mathbb{N}.$ The map $p$ must be such that $\sum_{j,k}p_{l,j,k}=n_l.$ Also we set $$m!=m_1!...m_d!,\ \  \overline{p}_l = (p_{l,1,1},...,p_{l,d,m_d}) \mbox{\ \ and \ \ } \overline{p}_{j,k}=(p_{l,j,k},...,p_{d,j,k}). $$   }
\end{enumerate}}
\end{defn}

For one variable, $d=1$, the conditions on the
coefficients are, respectively, as follows:
$$s_n(a+b)= \sum_{k=0}^{n}{n \choose k}s_k(a)
s_{n-k}(b),$$
$$s_n(a+b)=\sum_{i_1 + ... + i_k = n}\frac{1}{k!} {n \choose {i_1,\dots,i_k}} s_k(a)s_{i_1}(b)\dots s_{i_k}(b).$$ In this case these coefficient are called
binomial coefficients and multinomial coefficients, respectively
\cite{GCRota}.  Given $s={\des \sum_{n\in
\mathbb{N}^{d}} s_n
\frac{x^{n}}{n!}\in R[[x_1,\dots,x_d]]}$  we set $s^{0}=1$, and for $a\in \mathbb{N}$ we set $s^{a+1}=s^{a}s,$ and $s^{a}={\des \sum_{n\in
\mathbb{N}^{d}} s_n(a) \frac{x^{n}}{n!}}$. Also for $S\in C^{\mathbb{B}^{d}}$ we set $S^{0}=1$ and
$S^{a+1}=S^{a}S$.

\begin{thm}{\em
\begin{enumerate}
\item{ The sequence $\{s_n\}_{n\in\mathbb{N}^{d}}$ defined as above is polybinomial sequence.}
\item{ Assume that there exists $S\in C^{\mathbb{B}^{d}}$ such that
$|S|=s,$ then \\  $${\des s_n(a)=\left|
\bigoplus_{x_1\sqcup\dots \sqcup x_a=x} \bigotimes_{i=1}^{a} S(x_i,f|_{x_i})\right|}.$$ }
\item{Suppose that $S\in C^{\mathbb{B}^{d}}$ is such that $S=1 - F$ where $F \in C^{\mathbb{B}^{d}}_+$, then the
species $S^{-1}$ sending the empty set into $1$ and a non-empty set into
$$S^{-1}(x,f)= \bigoplus_{n \in \mathbb{N}}  \bigoplus_{x_1\sqcup\dots \sqcup x_n=x} \bigotimes_{i=1}^{n} F(x_i,f|_{x_i}) $$
is such that $|S^{-1}||S|=1=|S||S^{-1}|.$
}
\end{enumerate}}
\end{thm}

\bigskip

We define categories of directed trees $T$, $T_i(x,f),$ and
$T^{a}_i(x,f)$ pretty much as we did in the planar case, but now
we omit the planar condition. Let $s=(s_1,\dots,s_d)\in R[[
x_1,\dots,x_d]]^{d}$ be such that $s_i(0)=0$ and $\partial_j x_i = \delta_{i,j}$.
Set $s^{<0>}=(x_1,...,x_d)$, $s^{a+1}=s^{a}\circ s,$ and $$s^{\langle a
\rangle}_i={\des
\sum_{n\in\mathbb{N}^{d}} s_{n,i}(a)  \frac{x^{n}}{n!}}.$$
If $S\in C^{[d] \times\mathbb{B}^{d}}_+$ then we set $S^{0}=(X_1,...,X_d)$ and $S^{\langle
a\rangle}=S\circ S\circ\dots \circ S$.

\begin{thm}{\em
\begin{enumerate}
\item{ The sequence $\{ s_{n,i}\}$ defined above
is a  polymultinomial sequence. }
\item{Suppose that $S\in C^{[d] \times \mathbb{B}^{d}}_+$ satisfies $|S_i|=s_i,$ and let
$(x,f)\in \mathbb{B}^{d}$ be such that $f^{-1}(i)=n_i$. Then
$${\des s_{n,i}(a)=\left|\bigoplus_{\gamma\in \underline{T_i^{a}(x)}}
\left(\bigotimes_{v\in V\backslash x}S_{l(v)}(t^{-1}(v), l)\right)_{{\rm Aut}(\gamma)}   \right|.}$$}
\item{If $S\in C^{[d] \times \mathbb{B}^{d}}_+$ is such that $S_i=X_i-F_i$ with $F_{i}(x,f)=0$ for $|x|\leq 1$,
then the $C$-species $S^{<-1>}$ given by $${\des S^{<-1>}(x,f)=\left|\bigoplus_{\gamma\in \underline{T_i(x,f)}}
\left(\bigotimes_{v\in V_{\gamma}\backslash x}F_{l(v)}(t^{-1}(v), l)\right)_{{\rm Aut}(\gamma)}   \right|}$$
is such that $$|S|\circ |S^{<-1>}|=(X_1,...,X_d)= |S^{<-1>}|\circ |S| .$$}
\end{enumerate}}
\end{thm}
In the statement above the notation  $$ \left(\bigotimes_{v\in V_{\gamma}\backslash
x}S_{l(v)}(t^{-1}(v), l)\right)_{{\rm Aut}(\gamma)} $$ represents the
colimit, which we assume that exists, of the functor that sends
$\gamma\in T_i^a(x)$ into
$$\bigotimes_{v\in V_{\gamma}\backslash x}S_{l(v)}(t^{-1}(v),l).$$

\section{Categorification of Feynman integrals}\label{CFI}

In order to find an appropriated categorification of finite
dimensional Feynman integrals we need to generalize the notion of
species to the context of $G$-$C$ species; here $G$ is a
semisimple groupoid such that for $x\in G$ the cardinality of
$G(x,x)$ is finite, and $C$ is a symmetric distributive category
provided with a $R$-valuation. We call a groupoid semisimple if:
\begin{enumerate}
\item{It is provided with a bifunctor $\oplus: G \times G \longrightarrow G $, turning $G$ into a symmetric
monoidal category. We choose an unit object and denote it by $0$.}
\item{Every object in $G$ is isomorphic to a finite sum, unique up to reordering, of simple
objects. An object $x\in G$ is simple if $x=x_1\oplus x_2$
implies that either $x_1$ is isomorphic to $x$ and $x_2$ is isomorphic to $0$, or $x_1$ is isomorphic to $0$
and $x_2$ is isomorphic to $x$. }
\end{enumerate}

Let $G$ be a semisimple groupoid and $g_1,\dots,g_n,\dots$ a countable family
of formal variables.
A $R[[g_1,\dots,g_n,\dots]]$-weight on $G$ is a map
$\omega:G\longrightarrow R[[g_1,\dots,g_n,\dots]]$ satisfying
for  $x,y\in G$ the following conditions:
$$\omega(x)=\omega(y) \ \mbox{if} \ x  \ \mbox{is  \ isomorphic  \ to}  \ y,
 \ \omega(x\oplus y)=\omega (x)\omega(y),     \mbox{\ and\ }
\omega(0)=1.$$

We define the category of $G$-$C$ species to be the category $C^G$
of functors from $G$ to $C$. With these notation one defines a
map $|\ \ |:C^{G}\longrightarrow R[[g_1,\dots,g_n,\dots]]$ by
$$|F|=\sum_{x\in \underline{G}}
|F(x)|\frac{\omega(x)}{|G(x,x)|}.$$

Recall that a matching on a set $x$ is a partition of $x$ with blocks of cardinality two. Let $M:\mathbb{B}
\longrightarrow \mathbb{B}$ be the species that sends  $x$ to the set $M(x)$ of all matchings of $x$.
For applications to Feynman integrals we consider the groupoid of graphs $\graphs$. An  object $\gamma \in \graphs$
is a triple $(F, V,E)$ such that:
\begin{itemize}
\item{$F$ is a finite set whose elements are called flags.}
\item{$V$ is a partition of $F$; blocks of $V$ are the vertices of $\gamma$.}
\item{$E$ is a matching on $F$; blocks of $E$ are the edges of $\gamma$.}
\end{itemize}

\begin{figure}[h!]
\begin{center}
\includegraphics[width=6in]{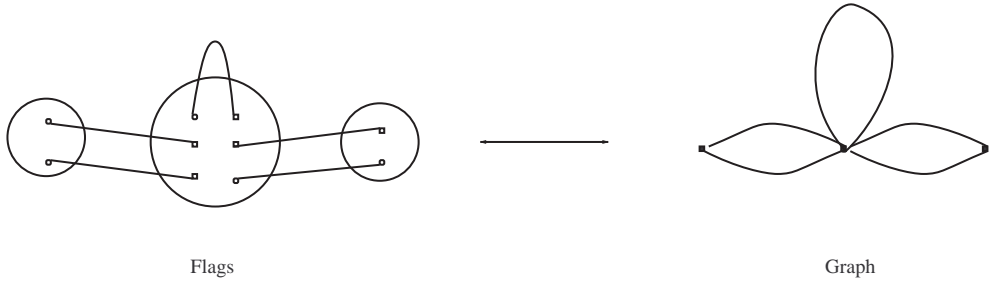}
\caption{Graphical representation of graphs. \label{fig:flags}}
\end{center}
\end{figure}

Figure \ref{fig:flags} shows how objects in $\graphs$ are
represented in pictures. A morphism
$\varphi:\gamma_1\longrightarrow
\gamma_2$ in $\graphs$ is a bijection $\varphi:F_{1} \longrightarrow F_{2}$ such that
$\varphi(V_{1})=V_{2}$ and $\varphi(E_{1})=E_{2}$. \\

On $C^{\graphs},$ the category of functors from $\graphs$ into $C$, we define the sum of functors as usual
$$(F+ G)(\gamma)=F(\gamma)\oplus G(\gamma).$$ The product functor on
$C^{\graphs}$ is given by
$$(F
G)(\gamma)={\des\bigoplus_{\gamma_1\sqcup\gamma_2=\gamma}F(\gamma_1)\otimes
G(\gamma_2)}$$ where $\gamma_1\sqcup\gamma_2$ denotes the disjoint
union of graphs. These definitions turn
$(C^{\graphs},+ , .)$ into a distributive category. Moreover each $R[[g_1,\dots,g_n,\dots]]$-weight on $\graphs$
induces a valuation map on $(C^{\graphs}, +, .)$.

\begin{prop}{\em The map  $|\ \ |:C^{\graphs}\longrightarrow
R[[g_1,\dots,g_n,\dots]]$ given by
$${\des|F|=\sum_{\gamma\in \underline{\graphs}}
|F(\gamma)|\frac{\omega(\gamma)}{|\graphs (\gamma,\gamma)|}}$$
defines a valuation on $C^{\graphs}$.}
\end{prop}

Notice that a vertex $v$ of a graph is a subset of the set
of flags, thus it makes sense to compute its cardinality  $|v|$. We shall use the following type of
$R[[g_1,\dots,g_n,\dots]]$-weight on $\graphs$:
$$\omega(\gamma)=\prod_{v\in V} g_{|v|}\ \ \mbox{for}\ \ \gamma\in \graphs.$$
Below we shall also need  the following map:
\[
\begin{array}{cccc}
  \pert: & \mathbb{C}[[x_1,\dots,x_d]]& \rightarrow & \mathbb{C}[[x_1,\dots,x_d,g_0,g_1,\dots,g_k,\dots]] \\
   & {\des\sum_{a\in\mathbb{N}^{d}}f_a x^{a}}& \mapsto &{\des \sum_{a\in\mathbb{N}^{d}}f_a x^{a}\frac{g_{|a|}}{|a|!}}  \\
\end{array}
\]
where for $a\in\mathbb{N}^{d}$ we set $|a|=\sum_{i\in [d]} a_i.$

A fundamental property of the  Gaussian measure is that it has a
clear combinatorial meaning; in contrast, a similar understanding
for the Lebesgue measure is lacking. The combinatorial meaning of
Gaussian integrals may be summarized in the remarkable identity:
\begin{equation*}\label{match}
\frac{1}{\sqrt{2\pi}}{\displaystyle\int_{-\infty}^{\infty}
x^{n}e^{-\frac{x^{2}}{2}}dx=|M[n]|}.
\end{equation*}
For example, see Figure \ref{match}, we have that
$$\frac{1}{\sqrt{2\pi}}{\displaystyle \int_{-\infty}^{\infty}
x^{4}e^{-\frac{x^{2}}{2}}dx=|M[4]|}=3.$$

\begin{figure}[h!]
\begin{center}
\includegraphics[width=11cm]{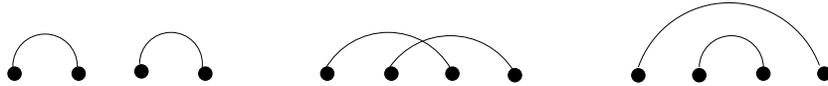}
\caption{Matchings on
$[4]$. \label{fig:match}}
\end{center}
\end{figure}

Using polarization and diagonalization, the formula above implies
the following identity for any positive definite symmetric
$n^{2}$\ real matrix $a$ and any $(x,f) \in
\mathbb{B}^{d}$:
\begin{equation*}
{\des \frac{(2\pi)^{-\frac{d}{2}}}{\sqrt{ {\rm det} a}}
\int_{\mathbb{R}^{d}} \prod_{i\in x} x_{f(i)}e^{-\sum a_{ij}^{-1} x_ix_j}dx_1...dx_d=
\sum_{\sigma\in M(x)}\prod_{m\in \sigma} a_{f(m)}},
\end{equation*}
where for $m=\{i,j\} \in \sigma$  we set $f(m)=\{f(i),f(j)\}$.\\

Assume we are given $A_{ij} \in C,$ for $i,j\in [d],$ such that
$A_{ij}\simeq A_{ji}$ and $|A_{ij}|=a_{ij}$. Given $(x,f)\in
\mathbb{B}^{n}$ we define $\graphs(x,f)$ to be the full subcategory
of $\graphs$ whose objects are graphs such that:
$x\subset F$ and $\{i\}\in V$ for all $i\in x$. Thus $\graphs(x,f)$ denotes
the category of graphs that include $x$ as a subset of the vertices of cardinality $1$.

\begin{defn}{\em Fix  $(x,f)\in
\mathbb{B}^{n}$. The Feynman functor $\mathfrak{F}: C^{\mathbb{B}^{d}}\longrightarrow C^{\graphs(x,f)}$
sends  $S\in C^{\mathbb{B}^{d}}$ to $\mathfrak{F}S\in C^{\graphs(x,f)}$
given by
$$\mathfrak{F}S(\gamma)={\des \bigoplus_{\widehat{f}}\left( \bigotimes_{v\in V }S(v,\widehat{f}|_v)
\bigotimes_{e\in E} A_{f(e)} \right) }\ \ \ \mbox{for }\ \ \gamma\in \graphs(x,f),$$
where the sum runs over the extensions $\widehat{f}:F\longrightarrow [d]$ of
$f:x \longrightarrow [d]$.}
\end{defn}

\begin{thm}{\em Fix $(x,f)\in
\mathbb{B}^{n}$ and let $S\in C^{\mathbb{B}^{d}}$ be such that $S(y,g)=0$ if $|y|\leq
2$. The following identity holds  $${\des
\frac{(2\pi)^{-\frac{d}{2}}}{\sqrt{ {\rm det} a}}\int_{\mathbb{R}^d}
e^{-a^{-1}_{ij}x_ix_j+\pert|S|}\prod_{i\in x}
x_{f(i)} dx}=|\mathfrak{F}S|.$$}
\end{thm}
\begin{proof}
Let the valuation of $S$ be given by
$$|S|=\sum_{t \in \mathbb{N}^d}s_t\frac{x^t}{t!} .$$

Then we have that
\begin{eqnarray*}
\mbox{}& \mbox{}& {\des\frac{(2\pi)^{-\frac{d}{2}}}{\sqrt{ {\rm det} a}}\int_{\mathbb{R}^d}
e^{-a^{-1}_{ij}x_ix_j+\sum_{|t|\geq 3}
s_t\frac{x^{t}}{|t|!}g_{|t|} }\prod_{i\in x}
x_{f(i)}dx_1...dx_d}=\nonumber\\
\mbox{}& \mbox{}& {\des\frac{(2\pi)^{-\frac{d}{2}}}{\sqrt{ {\rm det} a}}\int_{\mathbb{R}^d}
e^{-a^{-1}_{ij}x_ix_j} \prod_{|t|\geq 3}
\left(\sum_{k_t=0}^{\infty}
\frac{(s_t \frac{x^{t}}{|t|!} g_{|t|})^{k_t}}{k_{t}!}\right)
\prod_{i\in x}
x_{f(i)} dx_1...dx_d=}\nonumber\\
\mbox{}& \mbox{}& {\des\frac{(2\pi)^{-\frac{d}{2}}}{\sqrt{ {\rm det} a}}\int_{\mathbb{R}^d}
e^{-a^{-1}_{ij}x_ix_j} \prod_{|t|\geq 3}
\left(\sum_{k_t=0}^{\infty}
\frac{s_t^{k_t}  x^{k_t t} g_{|t|}^{k_t}}{k_{t}! |t|!^{k_t}}\right) \prod_{i\in x}
x_{f(i)} dx_1...dx_d}.
\end{eqnarray*}
Let $k$ be a map $k:\mathbb{N}^{d}\longrightarrow \mathbb{N},$ such $k_t=0$ for almost all $t \in \mathbb{N}^{d}$.
The formula above is equal to
\begin{equation*}\label{e8}
{\des\frac{(2\pi)^{-\frac{d}{2}}}{\sqrt{ {\rm det} a}}\sum_{k
}\left(
\prod_{|t|\geq 3}
\frac{s_t^{k_t}   g_{|t|}^{k_t}}{k_{t}!|t!|^{k_t}}\right)\int_{\mathbb{R}^d} e^{-a^{-1}_{ij}x_ix_j}  x^{\sum_{t} k_t
 t}\prod_{i\in x} x_{f(i)}dx_1...dx_d}.
\end{equation*} Let $z=y\sqcup x$ and
$(y,g)\sqcup(x,f)$ be such that $(y,g)$ is any colored set with
${\des g^{-1}(i)=\sum_{t=3}^{\infty} k_t t_i}$. Then the previous formula
becomes

$${\des \sum_{k}\sum_{\sigma\in M(z)}\left(
\prod_{|a|\geq 3}
\frac{s_t^{k_t}   g_{|t|}^{k_t}}{k_{t}!|t!|^{k_t}}\right)\prod_{m\in \sigma}
a_{f(m)}}={\des\sum_{\gamma\in \underline{\graphs}}|\mathfrak{F}S(\gamma)|
\frac{\omega(\gamma)}{\graphs(\gamma,\gamma)}=|\mathfrak{F}S|}.$$

\end{proof}
In the computations above we assumed that we could interchange
infinite sums and integrals; that is the formal step in the
definition of Feynman integrals. The formalism introduce in this
section will be further developed in \cite{cas1, cas2, RDEP1}.

\section{Categorification of deformation quantization}\label{kon}

In this section we assume that the reader is familiar with the
notations  and results from \cite{kont}. A Poisson manifold  is a
pair $(M,\{
\  , \  \})$ where $M$ is a $d$-dimensional smooth manifold provided
with a bracket $\{ \ , \  \}: C^{\infty}(M) \otimes C^{\infty}(M)
\longrightarrow C^{\infty}(M)$ satisfying for all $f,g,h\in C^{\infty}(M)$ the following identities:
\begin{enumerate}
\item{$\{f,g\}=-\{g,f\}$.}
\item{$\{f,gh\}=\{f,g\}h+g\{f,h\}$.}
\item{$\{f,\{g,h\}\}=\{\{f,g\},h\}+\{g,\{f,h\}\}$.}
\item{$\{ \ \  , \ \}$ is a local bidifferential operator.}
\end{enumerate}

The axioms above imply that the bracket can be written in local
coordinates as
$$\{f,g\}(x)=\sum_{i,j}\alpha^{ij}\partial_i f\partial_j g$$ where $\alpha^{ij}$ is an
antisymmetric $m^{2}$ matrix with entries in $C^{\infty}(M)$. The
bivector $$\alpha=\sum_{i,j}\alpha^{ij}\partial_i\otimes
\partial_j\in \Gamma(M,\bigwedge^2 TM)$$ is called the Poisson
bivector associated with the Poisson manifold $(M,\{ \ ,\  \})$.

\begin{defn}{\em Let $M$ be a Poisson manifold. A formal deformation of  $C^{\infty}(M)$   is a star product $\star:
C^{\infty}(M) [[\hbar]]\otimes_{\mathbb{R}[[\hbar]]}
C^{\infty}(M)[[\hbar]]
\longrightarrow C^{\infty}(M)[[\hbar]]$ such that:
\begin{enumerate}
\item{$\star$ is associative.}
\item{ $f\star
g=\displaystyle{\sum_{n=0}^{\infty}B_n(f,g)\hbar^{n}}$, where
$B_n(\ , \ )$ are bidifferential operators.} \item{$f\star
g=fg+\frac{1}{2}\{f,g\}\hbar+O(\hbar^{2})$, where $O(\hbar^{2})$
are terms of order $\hbar^{2}$.}
\end{enumerate}}
\end{defn}

Kontsevich \cite{kont} constructed a canonical $\star$-product for
any finite dimensional Poisson manifold. For the manifold
$(\mathbb{R}^{d},\alpha)$ with Poisson bivector $\alpha$ the
$\star$-product is given by the formula
\begin{equation*}\label{starp}
f\star g =\sum_{n=0}^{\infty}\frac{\hbar^{n}}{n!}\sum_{\gamma \in
\underline{G_{n,2}}}
\omega_\gamma B_{\gamma,\alpha}(f,g),
\end{equation*}
where $G_{n,2}$ is the category of admissible graphs and $\omega_{\gamma}$ are some
constants which are independent of $M$ and $\alpha$. Let us proceed to define in details
the category of admissible graphs.

\begin{defn}{\em For $k,n\in\mathbb{N}\times\mathbb{N}_{\geq 2}$ we let $G_{k,n}$
be the full subcategory of $\digraphs$ whose objects, called
admissible graphs of type $(k,n)$, are directed graphs $\gamma$
such that:

\begin{enumerate}
\item{$V_\gamma=V_\gamma^{1}\sqcup V_\gamma^{2}$ where $V_\gamma^{1}$ and $V_\gamma^{2}$ are totally ordered sets
with $|V_\gamma^{1}|=k$, $|V_\gamma^{2}|=n$.}
\item{$E_{\gamma}=V_{\gamma}^{1}\times [2]$.}
\item{$t(e)\neq s(e)$ for $e\in E_{\gamma}$.}
\item{$s(v,i)=v$ for $(v,i)\in E_{\gamma}$.}
\end{enumerate}}
\end{defn}

Next we define a couple of partial orders
$\leq_{L}$ and $\leq_{R}$  on  $\bigsqcup_{k\geq 0}G_{k,n}.$
First we need some graph theoretical notions.
Let $\gamma_1$ and
$\gamma_2$ be directed graphs. We say that $\gamma_{1}$ is
included in $\gamma_{2}$ and write $\gamma_{1}\subset
\gamma_{2}$, if $V_{\gamma_1}\subset V_{\gamma_2}$,
$E_{\gamma_1}\subset E_{\gamma_2}$ and
$(s_1,t_1)=(s_2|_{E_{\gamma_1}},t_2|_{E_{\gamma_1}})$. If
$\gamma_{1}\subset \gamma_{2}$ then we define the graph $\gamma_2/
\gamma_1$ as follows:
\begin{enumerate}
\item{ $V_{\gamma_2/
\gamma_1}=(V_{\gamma_2}\backslash V_{\gamma_1})\sqcup\{\ast\}$.}
\item{$E_{\gamma_2/
\gamma_1}=E_{\gamma_2}\backslash (s_2^{-1}(V_{\gamma_{1}})\cap t_2^{-1}(V_{\gamma_1}))$.}
\item{ ${\des s_{\gamma_2/
\gamma_1}=s_2}$.}
\item{$t_{\gamma_2/
\gamma_1}$ is equal to $t_2$ on $E_{\gamma_2/\gamma_1}\backslash t_2^{-1}(V_{\gamma_1})$ and to $\{\ast\}$
on $E_{\gamma_2/\gamma_1}\cap t_2^{-1}(V_{\gamma_1})$. }
\end{enumerate}

Let $\gamma\in \bigsqcup_{k\geq 0}G_{k,n}$ and assume that
$V_{\gamma}^{2}=\{ i_1< i_2\dots < i_n\}$.  For $\gamma_1\in
\bigsqcup_{k=0}^{m}G_{k,n}$ we let $\gamma_1\leq_{L} \gamma$ if
and only if: $$ \gamma_1\subset \gamma, \ \  V_{\gamma_1}^{0}=\{i_1<i_2< \dots < i_s\} \mbox{ for some }  s\leq n, \mbox{ and }
\gamma/\gamma_1  \mbox{ is an admissible graph. }$$

\noindent Similarly  we let $\gamma_1
\leq_{R} \gamma$ if and only if
$$\gamma_1\subset \gamma,  \ \ V_{\gamma_1}^{0}=\{i_{s}<i_{s+1}< \dots < i_n\} \mbox{ for some }  s\leq n,
\mbox{ and } \gamma/\gamma_1 \mbox{ is an admissible graph.}$$

\begin{figure}[h!]
\begin{center}
\includegraphics[width=12cm]{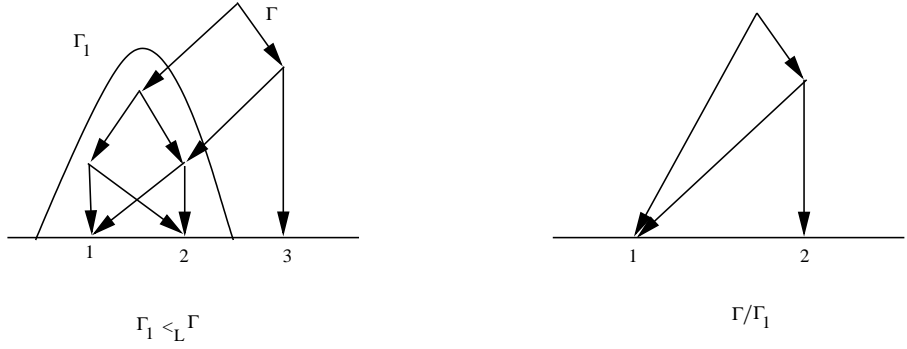}
\caption{Example of graphs such that $\gamma_1\leq_L \gamma$ and the associated quotient $\gamma/\gamma_1$. \label{fig:ss}}
\end{center}
\end{figure}

In order to categorify the Poisson manifold
$(\mathbb{R}^{d},\alpha)$ we need to find a distributive category with
a natural valuation on the  ring
$(R[[x_1,\dots,x_d,\hbar]],\star)$.  We make the following
assumptions:

\begin{enumerate}
\item{ A categorification $|\ \ |:C\longrightarrow R$  of $R$ is given.}
\item{For $1\leq i\neq j \leq d$ we are given $A^{ij}\in C^{\mathbb{B}^{d}}$ such that
$A^{ij}\simeq -A^{ji}$ and $|A^{ij}|=\alpha^{ij}$. }
\item{ We are given a functor $\Omega:\bigsqcup_{k\geq 0}G_{k,2}\longrightarrow C$ which sends $\gamma$ into
$\Omega_{\gamma}.$ We assume that there are natural isomorphisms:
\begin{equation*}\label{la}
{\des \sum_{\gamma_1\leq_{L}
\gamma}\Omega_{\gamma_1}\otimes\Omega_{\gamma/\gamma_1}\simeq
\sum_{\gamma_2 \leq_{R}\gamma}\Omega_{\gamma/\gamma_2}\otimes\Omega_{\gamma_2} }
\end{equation*}
for which Mac Lanes's pentagon axiom holds.}
\end{enumerate}

Recall that an object in $\mathbb{B}^{d+1}$ may be identified with
a triple $(x,f,y)$ where $x,y\in
\mathbb{B}$ and $f:x\longrightarrow [d]$.
\begin{thm}{\em
$(C^{\mathbb{B}^{d+1}}, +, \star)$ is a distributive category with
the $\star$-product on functors given by:
\begin{equation*}\label{estre}
F\star G(x,f,y)={\des
\bigoplus_{y_1\sqcup y_2\sqcup y_3=y}\
\bigoplus_{\gamma\in \underline{G_{|y_3|,2}}} \bigoplus_{I:E_{\gamma}\rightarrow [d]} \Omega_{\gamma}\otimes A_{\gamma,I,b}\otimes
B_{\gamma,I,b}}
\end{equation*}
where
\begin{enumerate}
\item{$I:E_{\gamma}\longrightarrow [d]$ and
$b:x\longrightarrow V_\gamma $ are maps.}
\item{$A_{\gamma,I,b}={\des
\bigotimes_{v\in V_{\gamma}}
A^{I(e_v^{1}) I(e_v^{2})}(b^{-1}(v)\sqcup E_v,f|_{b^{-1}(v)}\sqcup
I|_{E_{v}})}$.}
\item{$B_{\gamma,I,b}={\des
F(b^{-1}(\overline{1})\sqcup
E_{\overline{1}},f|_{b^{-1}(\overline{1})}\sqcup
I_{E_{\overline{1}}})\otimes G(b^{-1}(\overline{2})\sqcup
E_{\overline{2}},f|_{b^{-1}(\overline{2})}\sqcup
I_{E_{\overline{2}}})}$.}
\end{enumerate}}
\end{thm}
\begin{proof}
The key issue is that on the one hand we have that:
$${\des F\star(G\star H)(x,y,f)=\bigoplus \left(\bigoplus_{\gamma_1\leq_R \gamma} \Omega_{\gamma/\gamma_1}
\otimes\Omega_{\gamma_1}\right)\otimes A_{\gamma,I,b}\otimes B_{\gamma,I,b}}$$
where the sum runs over all $y_1\sqcup y_2\sqcup y_3= y$,
$\gamma\in \underline{G_{|y_3|,3}}$ and $I:E_{\gamma}\longrightarrow [m]$. On
the other hand we have that:
$${\des (F\star G)\star H(x,f,y)=\bigoplus \left(\bigoplus_{\gamma_1\leq_L \gamma} \Omega_{\gamma/\gamma_1}
\otimes\Omega_{\gamma_1}\right)\otimes A_{\gamma,I,b}\otimes B_{\gamma,I,b}}$$
where the sum runs over all $y_1\sqcup y_2\sqcup y_3= y$,
$\gamma\in \underline{G_{|y_3|,3}}$ and
$I:E_{\gamma}\longrightarrow [d]$. Above $A$,$B$ are defined as in
the statement of the Theorem \ref{tfn}.
\end{proof}

\begin{defn}{\em Let $(R[[x_1,\dots,x_d,\hbar]],\star)$ be the ring for formal power series in the variables
$x_1,\dots,x_d,\hbar$ with the star product $\star$ given by the Kontsevich's formula above with the constants $\omega_{\gamma}$ given by $\omega_{\gamma} = |\Omega_{\gamma}| \in R.$ }
\end{defn}

\begin{thm} {\em  The map $|\ \ |:(C^{\mathbb{B}^{d+1}}, +, \star)\longrightarrow
(R[[x_1,\dots,x_d,\hbar]],\star)$ given by $$|F|={\des
\sum_{(a,b)\in\mathbb{N}^{d}\times \mathbb{N}} |F([a],[b])|\frac{x^{a}}{a!} \frac{\hbar^{b}}{b!}}$$
defines a valuation  on $(C^{\mathbb{B}^{d+1}},+, \star)$.}
\end{thm}

\begin{proof} One checks that
$$F\star G={\des \sum_{n=0}^{\infty}\left(\sum_{\gamma\in \underline{G_{n,2}}}\Omega_{\gamma}
B_{\gamma}(F,G)\right)\frac{H^{n}}{n!}}$$ where $B_{\gamma}(F,G)$
is given by
$$\sum_{I:E_{\gamma}\longrightarrow [d]} \prod_{v\in V_{\gamma}^{0}}
\left( \prod_{e\in E_{\gamma}, t(e)=v} \partial_{I(e)}\right)A^{I(v,1)I(v,2)}
\left( \prod_{e\in E_{\gamma}, t(v)=1}\partial_{I(e)} \right)F
\left( \prod_{e\in E_{\gamma}, t(e)=v}\partial_{I(e)} \right)G.$$
Taking valuations and looking at page $5$ of \cite{kont} ones
obtains the desired result.
\end{proof}
\begin{figure}[h!]
\begin{center}
\includegraphics[width=9cm]{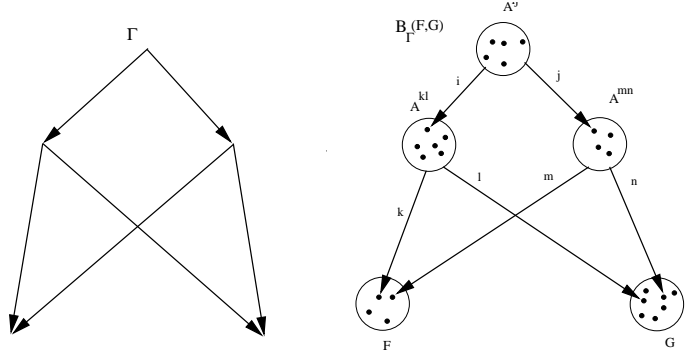}
\caption{ Example of a graph $\gamma\in \underline{G_{3,2}}$ and $B_{\gamma}(F,G)$. \label{fig:ss}}
\end{center}
\end{figure}

Finding a category $C$ with an appropriated family of objects $\Omega_{\gamma}$ is by no means an easy matter, fortunately Kontsevich's have shown that there are indeed examples \cite{kont}. We hope that the methods developed in this section may be of some use in order to find further examples.\\

Using  induction and the formula from Theorem \ref{estre}  one can show that:

\begin{thm}\label{tfn}{\em  Let $F_1\star F_2\star F_3\dots \star F_n \in C^{\mathbb{B}^{d+1}},$ then we have that:

$$F_1\star F_2\star F_3\dots \star F_n(x,f,h)\simeq
\bigoplus \bigoplus_{\gamma\in \underline{G_{|h_{n+1}|,n}}}\Omega_{\gamma}
\otimes A_{\gamma}\otimes B_{\gamma}$$ where the sum runs over the decompositions $h_1\sqcup\dots\sqcup
h_{n+1}=h$ of $h$ into $n+1$ disjoint blocks and
\begin{enumerate}
\item{${\des \bigoplus_{\gamma_1\leq_L\dots \leq_L\gamma_{n-2}\leq_L\gamma}\Omega_{\gamma_1}\otimes\bigotimes_{i=2}^{n-2}
\Omega_{\gamma_i/\gamma_{i-1}}\otimes \Omega_{\gamma /\gamma_{n-2}} }$.}
\item{$A_{\gamma}={\des
\bigotimes_{v\in V_{\gamma}}
A^{I(e_v^{1}) I(e_v^{2})}(b^{-1}(v)\sqcup E_v,f|_{b^{-1}(v)}\sqcup
I|_{E_{v}})}$.}
\item{ $B_{\gamma}={\des \bigotimes_{i=1}^{d}
F_i(b^{-1}(\overline{i})\sqcup
E_{\overline{i}},f|_{b^{-1}(\overline{i})}\sqcup
I_{E_{\overline{i}}})} $.}
\end{enumerate}}

\end{thm}

\section{Categorification of quantum phase space}\label{qspc}

In this section we consider the categorification of the quantum
phase space of a free particle with $n$-degrees of freedom. Before
developing the details of our approach we like to mention that
there are other attempts to try to understand quantum mechanics
using category theory, for example the reader may consult
\cite{mor, v}. Quantum phase space in this case is the deformation
quantization of the classical phase space, which may be identified
with the symplectic manifold $(\mathbb{R}^{2d}, \{ \ ,\  \})$ with
bracket:
$$\{x_i,x_j\}=\delta_{ij}, \ \ \ \{x_i,x_j\}=\{y_i,y_j\}=0$$ for
$x_1,\dots,x_d,y_1,\dots,y_d$ coordinates on $\mathbb{R}^{2n}$.
The ring of functions on the associated quantum phase space is
isomorphic to the formal Weyl algebra $W_{d}$ which we proceed to
introduce.

\begin{defn}{\em Let $R$ be a commutative ring. The Weyl algebra $W_{d}$ over $R$ is given by $$W_d=R\langle
\langle x_1,\dots,x_d,y_1,\dots,y_d \rangle
\rangle[[\hbar]]/I_d$$ where $I_{d}$ is the ideal generated for
$i,j\in [d]$ by
the following relations
$$\
[y_i,x_j]=\delta_{ij} \hbar, \ \ [x_i,x_j]=[y_i,y_j]=0.$$ }
\end{defn}

Objects in $\mathbb{B}^{2d+1}$ may be identified with triples
$(x,f,h)$ where $x, h\in \mathbb{B}$ and $f:x\longrightarrow
[2]\times [d]$. The color $(1,i)$ corresponds to the variable
$x_i$; the color $(2,i)$ corresponds to the variable $y_i$.

\begin{defn}\label{qps} {\em The distributive category  $(C^{\mathbb{B}^{2d+1}},+, \star)$
is such that the sum and product functors are given by: $$F\oplus G (x,f,h)=F(x,f,h)\oplus G(x,f,h),$$ $$F\star
G(x,f,h)=\bigoplus F (x_1 \sqcup h_3,f|_{x_1}\sqcup \{2\}\times
g,h_1)\otimes G(x_2 \sqcup h_3,f|_{x_2}\sqcup \{1\}\times g,h_2)$$
where the sum runs over all pairs $x_1, x_2$ and all triples $h_1,
h_2, h_3$ such that $$x_1\sqcup x_2=x, h_1\sqcup h_2\sqcup h_3= h
\mbox{ \ \ and \ \ } g:h_3 \longrightarrow [d].$$}
\end{defn}

\begin{figure}[h!]
\begin{center}
\includegraphics[width=14cm]{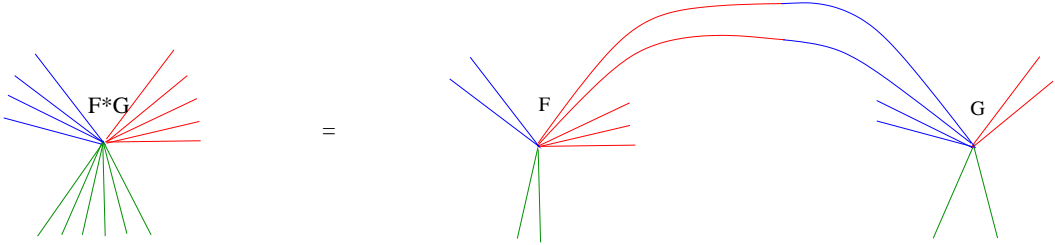}
\caption{Graphical interpretation of $F\star G$. \label{fig:starpro}}
\end{center}
\end{figure}

Figure \ref{fig:starpro} illustrates with an example the graphical
interpretation of the star product $F\star G$, where $F$ and $G$
are functors from $\mathbb{B}^{3}$ to $C$.\\

Our next result is a direct consequence of  the Proposition
\ref{propie} shown below.

\begin{thm} {\em The map $|\ \  |: C^{\mathbb{B}^{2d+1}}
\longrightarrow W_d$ given by
$$|F|={\displaystyle
\sum_{(a,b,c)\in \mathbb{N}^{d}\times\mathbb{N}^{d}\times \mathbb{N}}|F([a],[b],[c])|\frac{x^{a}}{a!}\frac{y^{b}}{b!}\frac{\hbar^{c}}{c!}}
$$
defines a valuation map on $(C^{\mathbb{B}^{2d+1}},+,
\star)$.}
\end{thm}

Using  Definition \ref{qps}  inductively ones obtains the
following result:

\begin{prop}\label{starm}{\em Given $F_1,\dots,F_m\in C^{\mathbb{B}^{2d+1}}$ we have that: $$F_1\star \dots \star
F_m(x,f,h)={\des\bigoplus
\bigotimes_{j=1}^{m}F_j(x_j
\bigsqcup_{i<j}h_{ij} \bigsqcup_{j<k}h_{jk}, f|_{x_j}\bigsqcup_{i<j}\{1\}\!\!\times\!\! g
\bigsqcup_{j<k}\{2\}\!\!\times\!\! g,h_{jj})}
$$ where the sum runs over the sets $x_i$, $h_{ij}$, and the maps $g$ such that:
$$\bigsqcup_{i=1}^{d} x_i=x, \ \
\bigsqcup_{i\leq j}h_{ij}=h \mbox{\ \  and \ \ }
g:\bigsqcup_{i<j}h_{ij}\rightarrow [d].$$}
\end{prop}

The result from the previous proposition may be rewritten as follows:
$$F_1\star \dots \star F_m(x,f,h)=\bigoplus_{\gamma\in G_{m}}
O_{\gamma},$$

where $G_m$ is the set of equivalence classes of  graphs $\gamma$
such that:
\begin{enumerate}
\item{$F_{\gamma}=x\sqcup h_1\sqcup [2]\times h_2$ where $h_1 \sqcup h_2=h$.}
\item{$g:h_2\longrightarrow [d]$ is any map. Below we use the
natural extension map $\widehat{g}:[2] \times h_2\longrightarrow
[2]\times [d]$.}
\item{$V_{\gamma}=\{b_1,\dots,b_d\}$ where $\bigsqcup_{i=1}^{d} b_i=F_{\gamma}$, $b_i\cap b_j=\emptyset$ and
if it happens that $s\in h_2$, $(2,s)\in b_i$, and $(1,s)\in b_j$,
then it must also happen that $i< j$.}
\item{$E_{\gamma}=\{(s,1),(s,2)|\ s\in h_2\}$.}
\end{enumerate}
We associate to each $\gamma\in G_m$ an object $O_{\gamma}\in C$
as follows:
$$O_{\gamma}={\des \bigotimes_{i=1}^{m} F_i(b_i\cap (x\sqcup [2] \times  h_2),
f|_{b_i\cap x}\sqcup \widehat{g}|_{b_i\cap ( [2]\times
h_2)},b_i\cap h_1)}.$$ With this notation it should be clear that
the formula above for the star product $F_1\star \dots \star F_m$
is just a reformulation of Proposition
\ref{starm}.\\

Next we consider the quantum analogue of the binomial sequences.

\begin{defn}\label{qps} {\em Let $m,n,l\in \mathbb{N}$. A sequence
$\{s_{m,n,l}\}$ where $s_{m,n,l}:\mathbb{N} \longrightarrow R$
such that
\begin{equation*}
s_{m,n,l}(r+t)= \sum_{l_1,l_2,l_3, m_1, m_2, n_1, n_2}  {l \choose
{l_1,l_2,l_3}}{m
\choose m_1}{n
\choose n_1} s_{m_1,n_1+s_3,l_1}(r)
s_{m_2+s_3,n_2,l_2}(t),
\end{equation*}
where $l_1 + l_2 + l_3=l,$ $m_1 + m_2 = m$ and $n_1 + n_2 = n$ is called a  quantum multinomial sequence.}
\end{defn}

Our next results describe a natural source of quantum multinomial
sequences and provide a categorical interpretation for such
sequences. Let
 $s={\des \sum_{a,b,c} s_{m,n,l}
\frac{x^{m}}{m!} \frac{y^{n}}{n!}\frac{\hbar^{l}}{l!}}\in W_d$, then we set $s^0=1,$ and for $r\in \mathbb{N}$ set
$s^{r+1}=s^{r} \star s$, and
$$s^{r}={\des \sum_{m,n,l} s_{m,n,l}(r) \
\frac{x^{m}}{m!} \frac{y^{n}}{n!}\frac{\hbar^{l}}{l!}}.$$
\begin{prop}{\em

\begin{enumerate}
\item{The sequence  $\{s_{m,n,l}\}$ defined  above is a
quantum multinomial sequence.}

\item{Assume that $S \in C^{\mathbb{B}^{3}}$ is such that $|S|=s,$ then
\begin{equation*}
s_{m,n,l}(r)=\left| {\des \bigoplus_{\Gamma\in G_n}
\bigotimes_{i=1}^{m} S(b_i\cap (x\sqcup [2] \times  h_2),
f|_{b_i\cap x}\sqcup
\widehat{g}|_{b_i\cap ( [2]\times h_2)},b_i\cap h_1)}\right|.
\end{equation*}}
\end{enumerate}}
\end{prop}

Let us consider a particular example of quantum multiplicative
sequence and provide a categorical interpretation for it.

\begin{defn}\label{llegando}{\em Let $n,a,b$ be integers such that $n-a-b$ is
zero or even. The quantum binomial is given  by $$
\Big[ \begin{array}{c}
  n \\
  a,b \\
\end{array}\Big]={\des \frac{n!(n-a-b-1)!!}{a!b!(n-a-b)!}}. $$}
\end{defn}

The reader should not confuse the integers $\Big[ \begin{array}{c} n \\
  a,b \\
\end{array}\Big]$ with the $q$-analogues of the binomial coefficients that are so often studied in the
literature. Below we need the singleton species $X$ and $Y$, they are define just as in the
case of commutative species.

\begin{prop}{\em
$ |(X+Y)^{n}([a],[b],[c])|= \Big[ \begin{array}{c}
  n \\
  a,b \\
\end{array}\Big] a!b!c!.$}

\end{prop}

\begin{proof}
According to Proposition \ref{starm} in order to construct
$(X+Y)^{n}([a],[b],[c])$, see Figure \ref{fig:star},  one should:
\begin{enumerate}
\item{Choose a partition of $[n]$ in $3$ blocks with cardinalities $a,b,$ and $2c$, respectively. There are
${\des{n\choose a,b,2c}}$ ways of doing this. Choose a linear
order on the first and second blocks. There are $a!b!$ such
possible orderings.}
\item{Choose a pairing on $[2c]$; there are $(2c-1)!!$ possible choices.}
\item{Choose a bijection between $[c]$ and the pairing selected in the previous step. There are $c!$ choices.}
\end{enumerate}
All together we see that
\begin{equation*}\label{bino}{\des |(X+Y)^{n}([a],[b],[c])|={n
\choose a,b,2c} a!b!(2c-1)!!}c! =\Big[ \begin{array}{c}
  n \\
  a,b \\
\end{array}\Big] a!b!c!.
\end{equation*}
\end{proof}
Using the formula above one obtains that:
$${\des (x+y)^{n}=|(X+Y)^{n}|=\sum_{a+b+2c=n}
 \frac{n!(2c-1)!!c!}{(2c)!} \frac{x^{a}}{a!}\frac{y^{b}}{b!}\frac{\hbar^{c}}{c!}=\sum_{a+b+2c=n}
\Big[ \begin{array}{c}
  n \\
  a,b \\
\end{array}\Big] x^{a}y^{b}\hbar^{c}}.$$

\begin{figure}[h!]
\begin{center}
\includegraphics[width=11cm]{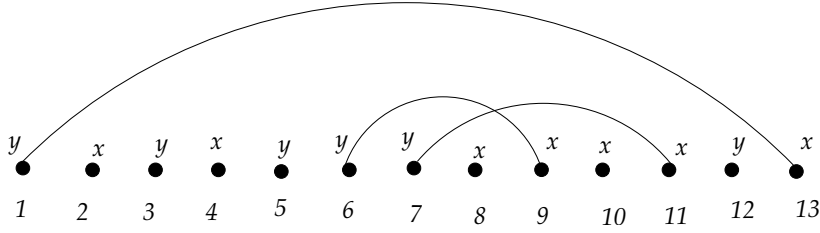}
\caption{Graph contributing to the computation of $(x+y)^{13}$ \label{fig:star}}
\end{center}
\end{figure}

The quantum binomial
coefficients satisfy a recursion relation which we describe next. It follows from identity  $(X+Y)^{n+1}=(X+Y)(X+Y)^{n}$
and identity $2.(b)$ below.

\begin{prop} {\em For  $n,a,b,c\in \mathbb{N}$ such that $a+b+2c=n+1$  the following
identity holds
\begin{equation*}\label{recur}\Big[ \begin{array}{c}
  n+1 \\
  a,b \\
\end{array}\Big]=\Big[ \begin{array}{c}
  n \\
  a-1,b \\
\end{array}\Big]
+ \Big[ \begin{array}{c}
  n \\
  a,b-1 \\
\end{array}\Big]+
(a+1)\Big[ \begin{array}{c}
  n \\
  a+1,b \\
\end{array}\Big].
\end{equation*}}

\end{prop}

\begin{prop}\label{propie}{\em The following formulae hold in $(C^{\mathbb{B}^{3}}, +,
\star)$:
\begin{enumerate}
\item{$\overbrace{X\star X\star \dots \star X}^{p} \star
\overbrace{Y\star Y\star \dots \star Y}^{q}=X^{p}Y^{q}$.}
\item{$Y X=XY+H$.}
\item{${\des\frac{Y^{n}}{n!}\star \frac{X^{m}}{m!}=\sum_{i=0}^{\min}\frac{X^{m-i}}{m-i!}
\frac{Y^{n-i}}{(n-i)!} \frac{H^{i}}{i!}}$, where $\min=\min\{n,,m\}$.}
\end{enumerate}}
\end{prop}

\begin{proof} $1.$ is obvious.
The second formula of Definition \ref{qps} implies that for all
$a,b,c\in\mathbb{B}$ the product $Y\star X(a,b,c)=1$ if $|a|=|b|=1$ and
$|c|=0$, or $|a|=|b|=0$ and $|c|=1$. Otherwise
$YX(a,b,c)=\emptyset$. This result exactly agree with
$(XY+H)(a,b,c)$ for  $a,b,c\in \mathbb{B}^{3}.$ \\

Recall that $\frac{X^{n}}{n!}(a,b,c)=1$ if $|a|=n$ and
$|b|=|c|=0$, and $0$ otherwise. Since
$$\frac{Y^{n}}{n!}\star\frac{X^{m}}{m!}(a,b,c)=\frac{Y^{n}}{n!}(b\cup c)\otimes \frac{X^{m}}{m!}(a\cup
c),$$ then ${\des
\frac{Y^{n}}{n!}\star\frac{X^{m}}{m!}(a,b,c)}=1$ if and only if
$|b|+|c|=n$ and $|a|+|c|=m$ and zero otherwise. Therefore ${\des
\frac{Y^{n}}{n!}\star\frac{X^{m}}{m!}(a,b,c)}=1$ if and only if $|b|=n-|c|$, $|a|=m-|c|$ and $|c|\leq\min \{n,m\}$, that is
$${\des\frac{Y^{n}}{n!}\star
\frac{X^{m}}{m!}=\sum_{i=0}^{\min}\frac{X^{m-i}}{m-i)!}
\frac{Y^{n-i}}{(n-i)!} \frac{H^{i}}{i!}}.$$

\end{proof}

Taking valuations in the Proposition above, setting $a=1$, $m=n$ for part $2$, and
$b=n$, $m=1$ for part $3$ ones gets the following corollary.

\begin{cor}\label{propie2}{\em The following identities hold in $W_1$:
\begin{enumerate}
\item{${\des\frac{y^{n}}{n!}\star \frac{x^{m}}{m!}=\sum_{i=0}^{\min}\frac{x^{m-i}}{m-i!}
\frac{y^{n-i}}{(n-i)!} \frac{\hbar^{i}}{i!}}$, where $\min=\min\{n,,m\}$.}
\item{$yx^{n}=x^{n}y+nx^{n-1}\hbar$.}
\item{$y^{n}x=xy^{n}+ny^{n-1}\hbar$.}
\end{enumerate}}
\end{cor}

Let us provide another application of our categorical
approach to Weyl algebras.

\begin{prop}{\em
\begin{enumerate}
\item{$[XY,X^{n}]=nX^{n}H.$}
\item{ $e^{Y} e^{X} =e^{X} e^{Y} e^{H}$.}
\end{enumerate}}
\end{prop}

\begin{proof}
\begin{enumerate}
\item{ We have
$$Y\star X^{n}(a,b,c)=\bigoplus Y(a_1,b_1\sqcup c_0,c_1)\otimes X^{n}(c_0\sqcup a_2, b_2,c_2)$$
where the sum runs over all set $a,b,c$ such that $a_1\sqcup
a_2=a$, $b_1\sqcup b_2=b$ and $c_0\sqcup c_1\sqcup c_2=c$.
\begin{itemize}
\item{$|b_1|=1, \ c_0=\emptyset$ implies that $c=\emptyset$. We
have the specie $X^{n}Y$.} \item{$|b_1|=\emptyset,  \ |c_0|=1, \ \
|a_2|=n-1$. In this case we have the specie $nYX^{n-1}H$.}
\end{itemize}
finally, we have $[XY,X^{n}]=nX^{n}H.$}
\item{$e^{Y}e^{X}(a,b,c)=\bigoplus e^{Y}(a, b\cup c)\otimes e^{X}(a\cup
c,b) =e^Y(b\cup c) e^{X} (a\cup c)=1_{\mathcal{C}} = e^X e^Y
e^{H}(a,b,c)$ for $a,b,c\in \mathbb{B}$. }
\end{enumerate}
\end{proof}

\begin{exmp}{\em
Consider the star product $Y^4\star X^3$ of the species $Y^4, X^3$. The posible
graphs arising in this case are shown in Figure \ref{fig:star}.
Thus we see that in the Weyl algebra the following identity
holds:
$$y^4\star x^3=|Y^4\star X^3|= x^3
y^4+12xy^{3}\hbar+36xy^{2}\hbar^{2}+24y\hbar^{3}.$$}
\end{exmp}

\begin{figure}[h!]
\begin{center}
\includegraphics[width=7cm]{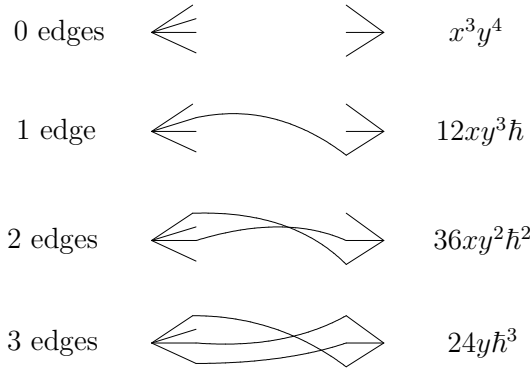}
\caption{ Graphs contributing to the computation of $y^4\star x^{3}$.\label{fig:star}}
\end{center}
\end{figure}

\section{Categorification of formal superspace}

In this closing section we study the categorification of the formal
supermanifold $\mathbb{R}^{d|m}$. Let $R$ be a commutative ring,
then by definition the ring of $R$-valued functions on
$\mathbb{R}^{d|m}$ is given by
$$R[[\mathbb{R}^{d|m}]]=R[[x_1,\dots,x_d]]\otimes
\bigwedge[\theta_1,\dots,\theta_m].$$ An element $f\in R[[\mathbb{R}^{d|m}]]$ can be uniquely written as
$f=\sum_{I\subset [m]}f_I \theta_{I}$ where $f_{I}\in
R[[x_1,\dots,x_d]]$, and for $I=\{i_1<\dots < i_k\}$ we set
$\theta_{I}=\theta_{i_1}\dots
\theta_{i_k}$. The product on $\mathbb{C}[[\mathbb{R}^{d|m}]]$ is given
by
\[
fg=\left(\des{\sum_{I\subset [m]}f_I \theta_{I}}\right)
\left(\des{\sum_{J\subset [m]}g_J \theta_{J}}\right)=
\des{\sum_{K\subset [m]}\left( \sum_{I\cup J=K}\sig(I,J) f_I g_J
\right)\theta_{K}},
\]
where $$ \ \ \theta_{I} \theta_{J}=\sig(I,J) \theta_{I\cup J} \mbox{\ \ \
and \ \ } \sig(I,J)=(-1)^{|\{ (i,j)\in I\times J|\  j<i
\}|}.$$ More generally if $I_1\sqcup I_2\sqcup\dots \sqcup I_p=I$
then we define
$$\sig(I_1, I_2, \dots,I_p)=\des (-1)^{|\{(i,j)|\ i\in I_k,
\  j\in I_l,\   1\leq k<l\leq p\  {\rm and} \ j<i|\}}.$$

\begin{defn}{\em
Given $d,m\in\mathbb{N}$ we let $\mathbb{B}^{d}\times F_m$ be the
category such that:
\begin{itemize}
\item{$\ob(\mathbb{B}^{d}\times F_m)=\{(x,f,I):
x\in \mathbb{B}^{d}, \ f:x\rightarrow [d], \ I\subset [m]\}$.}
\item{$\mathbb{B}^{d}\times F_m((x,f,I),(y,g,J))=\{ a: x \rightarrow
y\ | \ \alpha \mbox{\ is a bijection and \ } ga=f \}$ if $I= J$;
otherwise it is the empty set.}
\end{itemize}}
\end{defn}

\begin{defn}{\em
Let $C$ be a symmetric distributive category. The category $
C^{\mathbb{B}^{d}\times F_m}$ of functors from
$\mathbb{B}^{d}\times F_m$ to $C$ is called the category of
superspecies of type $d|m$.}
\end{defn}

Figure \ref{fig:ss} shows the graphical representation of the action of a
superspecies in $C^{\mathbb{B}^{5}\times F_8},$ where the standard
lines are bosons and the bold lines are fermions, i.e. represent commutative or anti-commutative variables.
\begin{figure}[h!]
\begin{center}
\includegraphics[width=7cm]{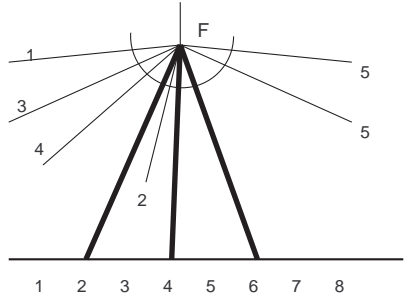}
\caption{ Graphical representation of a superspecies. \label{fig:ss}}
\end{center}
\end{figure}

\begin{defn}\label{propertss}{\em Let $C$ be a symmetric distributive
category and let $F,G\in C^{\mathbb{B}^{d}\times F_m}$. The
following formulae define the sum and product for superspecies:
$$(F + G)(x,f,I)=F(x,f,I)\oplus G(x,f,I),$$
$$F G(x,f,I)={\displaystyle
\bigoplus\sig(I_1,I_2)F(x_1,f|_{x_1},I_1)\otimes
G(x_2,f|_{x_2},I_2)}$$ where the sum runs over all $x_1\sqcup
x_2=x$ and $I_1\sqcup I_2=I$.}
\end{defn}

\begin{thm}{\em $(C^{\mathbb{B}^{d}\times F_m},+, .)$ is a symmetric distributive category, and  the map $$|\ \ |:
C^{\mathbb{B}^{d}\times F_m} \longrightarrow
R[[x_1,\dots,x_d]]\otimes
\bigwedge[\theta_1,\dots,\theta_m] $$  $$\mbox{given by} \ \ |F|=\sum_{(a,I)\in\mathbb{N}^{d}\times
F_m}|F([a],I)|\frac{x^{a}}{a!}\theta_{I}
\mbox{\ is \ a \ valuation \  map \ on \ }
C^{\mathbb{B}^{d}\times F_m}.$$}
\end{thm}

Next we define the super analogue of the binomial coefficients
\cite{GCRota}, and provided a combinatorial interpretation for them. We
also include a combinatorial interpretation for the multiplicative inverse of a
superspecies.

\begin{defn}{\em

Let $n\in\mathbb{N}^{d}$, $I\subset [m]$ and $a\in\mathbb{N}$.
A sequence $\{s_{n,I} \}$ where $s_{n,I}:\mathbb{N}
\longrightarrow R$ is called a super multiplicative sequence if
\begin{equation*}
s_{n,I}(a+b)={\des \sum_{j,A} \sig(A,A^{c}) {n\choose j}  s_{j,A}(a)
s_{n-j, I-A}(b)}
\end{equation*}
where $0 \leq j\leq n$ and $A\subset I$. }
\end{defn}

For $s={\des \sum_{n,I} s_{n,I}
\frac{x^{n}}{n!}}\theta_I\in R[[\mathbb{R}^{d|m}]]$ we set $s^{0}=1$, $s^{a+1}=s^{a}s,$ and
$$s^{a}={\des \sum_{n,I} s_{n,I}(a) \frac{x^{n}}{n!}}\theta_I.$$

\begin{prop}{\em \begin{enumerate}
\item{The sequence $\{s_{n,I}\}$ defined above is super multiplicative. }
\item{If $S \in C^{\mathbb{B}^{d}\times F_m}$ is such that $|S|=s,$ then
$$s_{n,I}(a)={\des|\bigoplus \sig(I_1,\dots,I_a) \bigotimes_{i=1}^{a} S(x_i,f|_{x_i},I_i)|}$$ where the
sum  runs over all partitions $x_{1}\sqcup\dots\sqcup x_a=x
\mbox{\ \  and \ \ }I_{1}\sqcup\dots\sqcup I_a=I.$}
\item{ If $S=1-F$ where $F\in C_+^{\mathbb{B}^{d}\times F_m}$ is such that $F(\emptyset)=0$,
then the superspecies
$$S^{-1}(x,f,I)={\des \bigoplus_{a=1}^{|x|+|I|} \bigoplus_{\sqcup_{1}^{a}x_i=x}\ \ \bigoplus_{\sqcup_{1}^{a}I_i=I}
\sig(I_1,\dots, I_a) \bigotimes_{i=1}^{a} F(x_i,f|_{x_i},I_i).
}$$ is such that $|S||S^{-1}|=1=|S^{-1}||S|.$ }
\end{enumerate}}
\end{prop}

The methods and techniques introduced in this work will gradually
find applications in a variety of settings. Applications of superspecies to the study of formal simple
supersymmetries will be developed in \cite{RDEP11}. For an
introduction to Lie algebras with a view towards categorification
the reader may consult \cite{dpgi}. One expects to find, along the
lines developed in this work, categorifications of several variants
of the Weyl algebra and their symmetric powers
\cite{DP0,DP2,RDEP3, DR}. It should  be possible to find a
categorical analogue of the perturbative methods developed by D\'iaz
and Leal \cite{dl} in order to obtain topological and geometrical
invariants from equivariant classical field theories.

\subsection*{Acknowledgments}
This work  owes  much to conversations of the first author with a
true teacher and friend Professor Gian-Carlo Rota. Part of this
work was done while the first author was visiting ICTP, Italy, and
Universidad de Sonora, M\'exico.

\noindent ragadiaz@gmail.com\\
\noindent Facultad de Administraci\'on, Universidad
del Rosario, Bogot\'a, Colombia \\

\noindent epariguan@javeriana.edu.co \\
Departamento de Matem\'aticas, Pontificia Universidad Javeriana,
Bogot\'a, Colombia

\end{document}